\documentclass[11pt]{article}
 \usepackage{amsfonts}
\usepackage{amsmath,amsfonts,amssymb,mathrsfs,cases,mathdots,colordvi,url,extarrows,graphicx}

\usepackage{relsize}
\usepackage{mathrsfs}
\usepackage{color}
\usepackage{indentfirst}
\usepackage{float}

\newtheorem{defn}{Definition}[section]

\newtheorem{lemma}{Lemma}[section]
\newtheorem{thm}{Theorem}[section]
\newtheorem{prop}{Proposition}[section]
\newtheorem{cor}{Corollary}[section]
\newtheorem{example}{Example}[section]
\newtheorem{remark}{Remark}[section]

\renewcommand{\Box}{\rule{2.2mm}{2.2mm}}

\def\beginproof{\par\noindent {\bf Proof.}\ \ }
\def\endproof{\hskip .5cm $\Box$ \vskip .5cm}

\def\beginproof{\par\noindent {\bf Proof.}\ \ }
\def\endproof{\hskip .5cm $\Box$ \vskip .5cm}

\begin{document}
\title{Directional quasi-/pseudo-normality conditions as sufficient conditions fro metric subregularity\thanks{The alphabetical order of the authors indicates the equal contribution to the paper.}}
\author{Kuang Bai\thanks{Department of Mathematics and Statistics, University of Victoria, Canada. The research of this author was partially supported by the China Scholarship Council. Email: kuangbai@uvic.ca.} \and Jane J. Ye \thanks{Corresponding author. Department of Mathematics and Statistics, University of Victoria, Canada. The research of this author was partially
supported by NSERC. Email: janeye@uvic.ca.}
\and Jin Zhang\thanks{Department of Mathematics, Southern University of Science and Technology, Shenzhen, P.R. China.
  This author's work is supported by NSFC (11601458, 11971220 and 11871269). Email: zhangj9@sustech.edu.cn.}} 
             \date{}
\maketitle
\begin{abstract}
	In this paper we study sufficient conditions for  metric subregularity of a set-valued map which is the sum of a single-valued  continuous map and a locally closed subset. First we derive a sufficient condition for metric subregularity which is weaker than the so-called first-order sufficient condition for metric subregularity (FOSCMS) by adding an extra sequential condition. Then
	we  introduce directional versions of quasi-normality and pseudo-normality which
	are stronger than the new {weak} sufficient condition for metric subregularity  but weaker than the  classical quasi-normality and pseudo-normality respectively.
	Moreover we introduce a nonsmooth version of the second-order sufficient condition for metric subregularity  and show that it  is a sufficient condition for  the new sufficient condition for metric {sub}regularity to hold. An example is used to illustrate that  directional pseduo-normality can be weaker than FOSCMS.  For the class of set-valued maps where the single-valued mapping is affine and the abstract set is the union of finitely many convex polyhedral sets, we show that  {pseudo-normality and hence} directional pseudo-normality holds automatically at each point of the graph. Finally we apply our results to complementarity and the Karush-Kuhn-Tucker systems.
	
	\vskip 10 true pt
	

	\vskip 10 true pt
	
	
\end{abstract}

\newpage
\section{Introduction}

In this paper, we study
stability analysis of the system of the form
\begin{equation} \label{GS}
P(x)\in \Lambda, \end{equation}
where $\mathscr{X},\mathscr{Y}$ are finite-dimensional Hilbert Spaces, $P:\mathscr{X}\rightarrow \mathscr{Y}$ is  continuous near the point of interest  and $\Lambda$ is a subset of $\mathscr{Y}$ which is closed near the point of interest. Throughout the paper, unless otherwise specified, we assume that $\mathscr{Y}$ is an $m$-dimensional Hilbert space with inner product $\langle \cdot, \cdot \rangle $ equipped with the orthogonal basis $\mathscr{E}=\{e_1,\dots, e_m\}$. Without loss of generality, throughout this paper for any $ y\in\mathscr Y$ we denote $\langle y,e_i\rangle$ by $y_i$, $i=1,\ldots,m$.

Since the set $\Lambda$ is not required to be convex,  the system represented by $P(x)\in \Lambda$ is very general and many systems can be formulated in this form. In particular, various variational inequalities/complementarity systems can be reformulated in this form. For example, consider the cone complementarity system defined as
$$ {\cal K} \ni \Phi(x) \perp \Psi(x)\in {\cal K} ,$$
where ${\cal K} $  is a convex cone in $\mathscr{Y}$, $\Phi, \Psi:\mathscr{X}\rightarrow \mathscr{Y}$, and $y \perp z$ means that  $\langle y, z \rangle =0$.  Then the cone complementarity system can be reformulated in the form $(\ref{GS})$ by defining  $P(x):=(\Phi(x), \Psi(x))$ and the complementarity set $$\Lambda:=\{ (y, z)\in \mathscr{Y} \times \mathscr{Y} | {\cal K} \ni y \perp z \in {\cal K} \}.$$
Note that although ${\cal K}$ is convex, the complementarity set  is not convex.

Denote by $G(x):=P(x)-\Lambda$,  a set-valued map induced by the system $P(x)\in \Lambda$. An important stability issue to study  is the  metric subregularity. We say that the set-valued map $G$ is metrically subregular at $(\bar{x}, 0) \in gph G$, where  $$gph G:=\{(x,y)| y\in G(x)\}$$ is the graph of $G$, if there {exist} $\kappa\geq 0$ and a neighborhood $U$ of $\bar x$ such that
\begin{equation*}\label{DefnMSR} d(x, G^{-1}(0))\leq \kappa d(P(x), \Lambda) \quad \forall x\in U,\end{equation*}
where $d(x,C)$ denotes the distance between a point $x$ and a set $C$ and $G^{-1}(y):=\{x| y\in G(x)\}$ denotes the inverse of $G$ at $y$.

The concept of metric subregularity was introduced by Ioffe \cite{Ioffe79} using the terminology ``regularity at a point.'' The terminology ``metric subregularity'' was suggested  by Dontchev and Rockafellar in \cite[Definition 3.1]{DonRock04}.
This property is also referred to as an error bound property since it enables us to estimate the distance from a point $x$ near $\bar x$ to the set of solutions to the system $(\ref{GS})$ by its residue $d(P(x),\Lambda)$,  which is much easier to deal with; see, e.g., \cite{P,WY01,WY02,WY03,FHKO,Penot} and  the references therein for related results and applications.
Metric subregularity is a weaker condition than the more familiar property of metric regularity which requires the existence of $\kappa \geq 0$ and $U, V$, neighborhoods of $\bar x, 0$, respectively, such that
$$ d(x, G^{-1}(y))\leq \kappa d(P(x), \Lambda) \quad \forall x\in U, y\in V,$$
and strong metric subregularity (see, e.g., \cite{DonRock}) which requires the existence of $\kappa \geq 0$ and $U$, a neighborhood of $\bar x$ such that
$$ \|x-\bar x\| \leq \kappa d(P(x), \Lambda) \quad \forall x\in U.$$

 It is well known (see e.g. \cite[Theorem 3.2]{DonRock04})  that the metric subregularity of a set-valued map  is equivalent to the calmness  of its inverse map, which means that  there exist $\kappa\geq 0$ and neighborhoods $U$ of $\bar x$ and $V$ of $0$ such that
$$G^{-1}(y)\cap U\subseteq G^{-1}(0) +\kappa \|y\| \mathbb{B} \quad \forall y \in V,$$
where $\|\cdot\|$ and ${\mathbb{B}}$ denote the norm and  the closed unit ball in $\mathscr{Y}$, respectively.
The concept of the calmness was first introduced by J. J. Ye and X. Y. Ye in  \cite[Definition 2.8]{YY} under a different name, ``pseudo upper-Lipschitz continuity,'' and the terminology of ``calmness'' was coined by Rockafellar and Wets in \cite{RW}. Note that the calmness property is part of the property required in   the notion of pseudo-Lipschitz continuity introduced by Klatte \cite{Klatte87}. 
As  suggested by the name ``pseudo upper-Lipschitz continuity,'' the concept of calmness is weaker than both the pseudo-Lipschitz continuity (or Aubin continuity) introduced by Aubin \cite{Aub} and the upper-Lipschitz continuity introduced by Robinson \cite{Robinson75}-\cite{Robinson81}. Analogous to the fact that a set-valued map is metrically  subregular  if and only if its inverse map is calm, it is well known that the metric regularity of a set-valued map is equivalent to the pseudo-Lipschitz continuity of its inverse map (see \cite[Theorem 1.49]{Aub2}).

Metric subregularity/calmness plays an important role in optimization. It serves as a constraint qualification and a sufficient condition for exact penalty; see e.g., \cite{Burke,Ioffe79,Ioffe79new,Ioffe00,KK,SW,Ye,YY}. As pointed out in
\cite{IO}, metric subregularity/calmness  is also an important tool in the  subdifferential calculus of nonsmooth analysis.  More recently, it has been discovered that it serves as a sufficient condition for linear convergence of certain numerical algorithms \cite{KK09,WYYZZ} and quadratic convergence of the Newton-type method \cite{FFH}.

Although the metric subregularity/calmness/error bound condition is very important, it is by no means easy to verify. For a long time,  there have been only two major {\it checkable} sufficient conditions:  one is derived by Robinson's multifunction theory and the other is by Mordukhovich's criteria. By Robinson's {multifunction} theory \cite{Robinson81}, if the linear constraint qualification (linear CQ) holds, i.e., $P(x)$ is affine and $\Lambda$ is the union of finitely many polyhedral convex sets, then the set-valued map $G(x)=P(x)-\Lambda$ must be a polyhedral multivalued function and so is its inverse map $G^{-1}$. Hence the set-valued map $G^{-1}$ must be upper Lipschitz and hence calm. {Recall that in optimization we call a multiplier abnormal if it is a multiplier corresponding to  an optimality system where the objective function vanishes.}  Assuming $P$ is continuously differentiable ($C^1$),  if  the no nonzero abnormal multiplier constraint qualification (NNAMCQ) holds, i.e., there is no nonzero  abnormal multiplier $\zeta$ such that
\begin{equation}\label{KKT}
0= \nabla P(\bar x)^* \zeta, \quad \zeta \in N_\Lambda(P(\bar x)) ,\end{equation}
where $N_\Lambda(\cdot)$ is the limiting normal cone, $\nabla P$ denotes the Fr\'{e}chet derivative  of $P$, and $^*$ denotes the {adjoint}, then the Mordukhovich's criteria for metric regularity (see, e.g., \cite[Theorem 9.40]{RW})  holds and so does metric subregularity. These two criteria are relatively strong since they are actually sufficient conditions for stronger stability concepts. And therefore there are many situations where these sufficient conditions do not hold but the systems are still metrically subregular. In general metric subregularity is weaker than NNAMCQ but for the case of differentiable convex inequality system, which is (\ref{GS}) with $P$ {convex and }differentiable and $\Lambda$ a nonnegative orthant, Li \cite{WuLi} has shown that all the following conditions are equivalent: metric subreguality, Abadie's constraint qualification, Slater condition and the  Mangasarian-Fromovitz constraint qualification (MFCQ) (which is equivalent to NNAMCQ in this case).

Over the last fifteen years or so,  some results  for  characterizing metric
subregularity/calmness for general set-valued maps have been obtained; see, e.g., \cite{HenJou,HenJouOutr,HenOutr01,HenOutr05,ZhengNg10}. Recently the concept of a directional limiting normal cone which is in general a smaller set than the limiting normal cone was introduced \cite{GM,Gfr13}. Based on the result for general set-valued maps in \cite{Gfr13},  Gfrerer and Klatte \cite[Corollary 1]{GKlatte16}  showed  that metric subregularity holds for system $(\ref{GS})$ at $\bar x$ under  the first-order  sufficient condition for metric subregularity (FOSCMS): assuming $P(x)$ is $C^1$,  if  for each nonzero direction $ u  $ satisfying $\nabla P(\bar x)u\in T_\Lambda(P(\bar x))$, there is no nonzero $\zeta$ such that
$$0= \nabla P(\bar x)^* \zeta, \quad \zeta \in N_\Lambda(P(\bar x); \nabla P(\bar x)u),$$
where $T_\Lambda(\cdot)$ and $N_\Lambda(y; d)$ are the tangent cone and the limiting normal cone at $y$ in direction $d$ (see Definition \ref{directionNormalC}). Moreover if $P(x)$ is strictly differentiable and twice directionally differentiable and $\Lambda$ is the union of finitely many polyhedral convex sets, it was shown in \cite[Theorem 4.3]{Gfr132} that  the metric subregularity holds at $(\bar x,0)$ under the following second-order sufficient condition for metric subregularity (SOSCMS):
 for each nonzero direction $ u  $ satisfying $\nabla P(\bar x)u\in T_\Lambda(P(\bar x)$,  there exists no $\zeta\neq0$ such that
	\begin{eqnarray*}
		\left \{  \begin{array}{l}
			0= \nabla P(\bar{x})^*\zeta, \quad  \zeta\in N_\Lambda(P(\bar{x});\nabla P(\bar{x})u),\\
 \langle \zeta,  P''(\bar{x};u) \rangle \geq 0
		\end{array}\right.
	\end{eqnarray*}
	where $P''(\bar{x};u)$ denotes the second-order derivative of $P(x)$ at $\bar{x}$ in the direction $u$. Some sufficient conditions for the metric subregularity/calmness/error bound condition  for special complementarity systems based on the FOSCMS have been obtained in \cite{GY,Yezhou17}.
	
Another direction in the effort of weakening the NNAMCQ is to add some extra conditions to $(\ref{KKT})$. 	In the case where $P$ is continuously differentiable at $\bar x$, we say that quasi-normality and pseudo-normality hold at  $\bar x$ if there exists no $\zeta\neq0$ such that $(\ref{KKT})$ holds and
$$ \exists (x^k, s^k, \zeta^k) \rightarrow (\bar{x}, P(\bar x),\zeta)~  \mbox{ s.t. }
 \zeta^k \in N_\Lambda(s^k)  \mbox{ and }  \zeta_i(P_i(x^k)-s_i^k) >0, \mbox{ if } \zeta_i \not =0,$$
$$ \exists (x^k, s^k, \zeta^k) \rightarrow (\bar{x}, P(\bar x),\zeta)~  \mbox{ s.t. }
 \zeta^k \in N_\Lambda (s^k)  \mbox{ and } \langle  \zeta,  P(x^k)-s^k\rangle >0 ,$$
 respectively. It is obvious that  pseudo-normality implies quasi-normality. For a system with equality and inequality constraints where all constraint functions are $C^{1+}$ which means that the gradients are locally Lipschitz, Minchenko and Tarakanov \cite[Theorem 2.1]{MinT}  showed that quasi-normality implies the existence of a local error bound or equivalently metric subregularity/calmness at $\bar x$. {In \cite[Theorem 5]{YZ}, this result is extended to systems with continuously differentiable equality constraint functions and subdifferentially regular inequality constraint functions and a  regular constraint set.}  Quasi-normality/pseudo-normality for the general system in the form $(\ref{GS})$ was introduced by
 Guo, Ye and Zhang  \cite[Definition 4.2]{Guo13} and proved to be a sufficient condition for error bound/metric subregularity/calmness in \cite[Theorem 5.2]{Guo13}  {under the Lipschitz continuity of $P$ and the closeness of the set $\Lambda$ only}.

The main purpose of this paper is to combine the two approaches of weakening the NNAMCQ , i.e., to replace the limiting normal cone by the directional normal cone as in  FOSCMS and SOSCMS and to add extra conditions as in quasi-/pseudo-normality and prove that our weaker sufficient conditions are still sufficient for verifying the  metric subregularity/calmness.

Our assumptions are very general. We only assume the continuity of the mapping $P(x)$. Indeed, it is {natural} to study the case where $P(x)$ is only continuous, since it will widen the range of applications of formation $(\ref{GS})$. For example, consider the recovery of an unknown vector $x \in \mathbb{R}^n$ (such
as a signal or an image) from noisy data $b \in \mathbb{R}^m$ by minimizing with respect to $x$ a regularized cost function
\begin{eqnarray}\label{nonsmoothDatafide}
F(x, b)= f(x, b) + \mu g(x),
\end{eqnarray}
where typically $f: \mathbb{R}^n \times \mathbb{R}^m \rightarrow \mathbb{R}$ is a data-fidelity term and
$g: \mathbb{R}^n \rightarrow \mathbb{R}$ is a nonsmooth regularization term, with $\mu > 0$ a parameter. One usual choice for the data-fidelity term is
$$f(x, b) = \sum_{i = 1}^m |a_i^T x  - b_i|^\rho$$
with $a_i\in \mathbb{R}^n$ and $\rho$ in the range $(0, \infty]$; see, e.g., \cite{app-nonlip1,app-nonlip3,app-nonlip2}. Apparently when $\rho$ takes a value in the interval $(1, 2)$, the optimality condition of minimizing function $(\ref{nonsmoothDatafide})$ with respect to $x$
 can be described by  $0\in \nabla_x f(x,b) +\partial g(x)$, where $\partial g(x)$ denotes a certain subdifferential of $g$ at $x$, which can be reformulated as $P(x)\in \Lambda$ where $P(x)$ is continuous. Therefore, thanks to the equivalence between calmness of different reformulations established in \cite[Proposition 3]{GY}, our results can be used to study the calmness of the optimality condition system of minimizing $(\ref{nonsmoothDatafide})$ without imposing an unnecessarily stronger condition.

We organize our paper as follows. Section 2 gives the preliminaries and preliminary results. In section 3, we propose the weak sufficient condition for metric subregularity and show that it is sufficient for metric subregularity. In Section 4, we propose the concepts of  directional quasi-/pseudo-normality and show that they are stronger than the new sufficient condition for metric subregularity. Moreover in this section it is shown that the SOSCMS implies pseudo normality. In section 5 we apply our results to complementarity systems and Karush-Kuhn-Tucker (KKT) systems.
	
\section{Preliminaries and preliminary results}
In this section, we gather some preliminaries on variational analysis and nonsmooth analysis that will be used in the following sections. We only give concise definitions and results that will be needed in this paper. For more detailed information on the subject, the reader is referred to Mordukhovich \cite{Aub2}, and Rockafellar and Wets \cite{RW}.

First, we give the definition of tangent cones and normal cones.
\begin{defn}
(tangent cones and normal cones; see, e.g. \cite[Definition 6.1]{RW}).
	Given a set $\Omega\subseteq\mathscr{Y}$ and a point $\bar{y}\in \Omega$, the tangent cone to $\Omega$ at $\bar{y}$ is defined as
	$$T_\Omega(\bar{y}):=\left \{d\in\mathscr{Y}|\exists t_k\downarrow0, d_k\rightarrow d\ \mbox{ s.t. } \bar{y}+t_kd_k\in \Omega\ \forall k\right \}.$$
	The derivable cone to $\Omega$ at $\bar{y}$ is defined as
	$$T_\Omega^i(\bar{y}):=\left \{d\in\mathscr{Y}|\forall  t_k\downarrow0, \exists d_k\rightarrow d\ \mbox{ s.t.  } {\bar{y}}+t_kd_k\in \Omega\  \forall k \right \}.$$
	A set $\Omega$ is said to be geometrically derivable  if the tangent cone coincides with the derivable cone at  each  point of $\Omega$, or equivalently
	if $\lim_{t\downarrow 0} t^{-1} d(\bar{y}+tu, \Omega)=0$.

	The regular normal cone and the limiting normal cone to $\Omega$ at $\bar{y}$ are defined as
	$$\widehat{N}_\Omega(\bar{y}):=\left \{\zeta \in  {\mathscr{Y}}\bigg| \limsup_{y\xrightarrow{\Omega}\bar{y}}\frac{\langle \zeta,y-\bar{y}\rangle}{\|y-\bar{y}\|}\leq0\right \},$$
	and $$ N_\Omega(\bar{y}):=\left \{\zeta\in {\mathscr{Y}}\bigg| \exists \ y_k\xrightarrow{\Omega}\bar{y},\ \zeta_k{\rightarrow}\zeta\ \text{such that}\ \zeta_k\in\widehat{N}_\Omega(y_k)\ \forall k\right \}$$
	respectively, where $y_k\xrightarrow{\Omega}\bar{y}$ means $y_k\rightarrow\bar{y}$ and for each $k$, $y_k\in \Omega$.
	\end{defn}
	Recently a directional version of limiting normal cones was introduced in \cite[Definition 2.3]{GM} and extended to   general Banach spaces in \cite{Gfr13}.
	\begin{defn}\mbox{(directional normal cones; see \cite[Definition 2]{Gfr13}).} \label{directionNormalC}
    Given a point $\bar y \in\mathscr{Y}$ and a direction $d\in\mathscr{Y}$, the limiting normal cone to $\Omega$ at $\bar{y}$ in direction $d$ is defined by
    $$N_\Omega(\bar{y};d):=\left \{\zeta \in \mathscr{Y}\bigg| \exists \ t_k\downarrow0, d_k\rightarrow d, \zeta_k\rightarrow\zeta  \mbox{ s.t. } \zeta_k\in \widehat{N}_\Omega(\bar{y}+t_kd_k)\ \forall k \right \}.$$
\end{defn}
From the definition, it is obvious that $N_{\Omega}(\bar y; d)=\emptyset$ if $d \not \in T_\Omega(\bar y)$ and $N_{\Omega}(\bar y;d)\subseteq N_\Omega(\bar y)$.

\begin{prop}\label{prop YZ}\cite[Proposition 3.3]{Yezhou17}  Let $\Omega:=\Omega_1\times\cdots\times \Omega_l$, where $\Omega_i
\subseteq \mathbb R^{n_i}$ are closed for $i=1,\ldots,l$ and $n=n_1+\dots+ n_l$.  Consider a point $\bar y=(\bar y_1,\ldots,\bar y_l)\in \Omega$ and a direction $d=(d_1,\ldots,d_l)\in\mathbb R^n$. Then
\begin{align*}
T_{\Omega}(\bar y)\subseteq T_{\Omega_1}(\bar y_1)\times\cdots\times T_{\Omega_l}(\bar y_l),\\
N_{\Omega}(\bar y; d)\subseteq N_{\Omega_1}(\bar y_1;d_1)\times\cdots\times N_{\Omega_l}(\bar y_l;d_l).
\end{align*}
	The equality holds if all except at most one of $\Omega_i$ for $i=1,\dots, l$ are directionally regular at $y_i$ in the sense of \cite[Definition 3.3]{Yezhou17}.
\end{prop}

We give the definition of some subdifferentials below.
\begin{defn}(subdifferentials; see, e.g., \cite{Aub2})
Let $f:\mathscr{X}\rightarrow [-\infty, +\infty]$ and $\bar x$ is a point where $f$ is finite. Then
\begin{itemize}
\item the Fr$\acute{e}$chet (regular) subdifferential of $f$ at $\bar{x}$ is the set
\begin{eqnarray*}
\widehat\partial f(\bar{x}):=\left\{\xi\in{\mathscr{X}}\bigg|\lim\inf_{h\rightarrow0}
\frac{f(\bar{x}+h)-f(\bar{x})-\langle\xi,h\rangle}{\|h\|}\geq0\right\};
\end{eqnarray*}
\item the limiting (Mordukhovich or basic) subdifferential of $f$ at $\bar{x}$ is the set
\begin{align*}
\partial f(\bar{x})&:=
\left \{\xi\in\mathscr{X}\big|\exists x_k\rightarrow \bar{x}, f(x_k)\rightarrow f(\bar x), \mbox{and}\ \xi_k\rightarrow\xi\ \mbox{with} \ \xi_k\in\widehat\partial f(x_k) \right\}  ;
\end{align*}
\end{itemize}
\end{defn}

Recently based on the concept of the directional limiting normal cone, the following directional version of  the limiting subdifferential was introduced in \cite{BGO}.
\begin{defn}(directional subdifferentials; see \cite{BGO})
Let $f:\mathcal X\rightarrow [-\infty, +\infty]$ and $\bar x$ be a point where $f$ is finite. Then the limiting subdifferential of $f$ at $\bar x$ in direction $(u, \zeta) \in \mathscr{X}\times \mathbb{R}$ is defined as
$$\partial f(\bar{x};(u,\zeta)):=
\left\{\xi\in\mathscr{X}\bigg|\exists t_k\downarrow0, u^k\rightarrow u, \zeta^k\rightarrow\zeta, \xi^k\rightarrow\xi,
f(\bar x)+t_k\zeta^k=f(\bar x+t_ku^k), \xi^k\in\widehat\partial f(\bar x+t_ku^k)\right \}.$$
\end{defn}
\begin{remark}
	Let $f$ be continuously differentiable at $\bar x$. Then $\partial f(\bar x;(u,\zeta))\neq\emptyset$ if and only if $\zeta=\nabla f(\bar x)u$, in which case
	$$\partial f(\bar x;(u,\zeta))=\partial f(\bar x)=\{\nabla f(\bar x)\}.$$
\end{remark}

\begin{defn}(graphical derivatives; see, e.g., \cite{DonRock})
For a set-valued map $G:\mathscr{X}\rightrightarrows\mathscr{Y}$ and a pair $(x,y)$ with $y\in G(x)$, the graphical derivative of $G$ at $x$ for $y$ is the set-valued map $DG(x|y):\mathscr{X}\rightrightarrows\mathscr{Y}$ whose graph is the tangent cone  to $gphG$ at $(x,y)$:
$$v\in DG(x|y)(u)\Leftrightarrow(u,v)\in T_{gphG}(x,y).$$
Thus, $v\in DG(x|y)(u)$ if and only if there exist sequences $u_k\rightarrow u,\ v_k\rightarrow v$ and $\tau_k\downarrow0$ such that $y+\tau_kv_k\in G(x+\tau_ku_k)$ for all $k$.
\end{defn}
For a single-valued mapping $P:\mathscr{X} \rightarrow\mathscr{Y}$, its graphical derivative at $x$ for $y=P(x)$ is
\begin{equation}
    D P(x)(u):=\left\{\xi\bigg|\exists t_k\downarrow0, u_k\rightarrow u \ \mbox{s.t.} \lim_{k\rightarrow+\infty}\frac{P(x+t_ku_k)-P(x)}{t_k}=\xi\right\}.
\end{equation}

Moreover if $P(x)$ is Hadamard directionally differentiable at $x$, then its graphical derivative is equal to the directional derivative: for any $u \in \mathscr{X}$,
$$D P(x)(u)=P'(x;u):=\lim_{t\downarrow 0, u'\rightarrow u}\frac{P(x+tu')-P(x)}{t}.$$

The following sum rule extends the sum rule in \cite[Proposition 4A.2]{DonRock} by allowing $P(x)$ to be only continuous.
\begin{prop}\label{prop gd}
Let $G(x):=-P(x)+\Lambda$ and $P(\bar x)\in \Lambda$, where $P(x):\mathscr{X}\rightarrow \mathscr{Y}$ is a continuous singled-valued map.

Then  either
\begin{equation}\label{lc}  DG(\bar{x}|0)(u)\subseteq -DP(\bar x)(u)+T_\Lambda(P(\bar{x}))\end{equation}
or there exists $\zeta\not =0$ such that
$$ \zeta \in DP(\bar x)(0)\cap T_\Lambda(P(\bar{x})).$$
If either $P(x)$ is Hadamard directionally differentiable  at $\bar x$ or $\Lambda$ is geometrically derivable, then  (\ref{lc})  holds as an equality.
\end{prop}
\beginproof
By  definition,  $v\in DG(\bar{x}|0)(u)$ if and only if $(u, v) \in T_{gph G} (\bar x,0)$.
It follows from  the definition of tangent cone that  there exist sequences $(u_k, v_k) \rightarrow (u,v)$ and $\tau_k\downarrow 0$ such that $(\bar x,0)+\tau_k(u_k,v_k) \in gph G$, which means that there exists $ s_k\in \Lambda$ such that $\tau_kv_k=-P(\bar{x}+\tau_ku_k)+s_k.$

 $Case$ (i) ($ \{ \frac{P(\bar{x}+\tau_ku_k)-P(\bar{x})}{\tau_k} \}$ is bounded.) Then without loss of generality we may assume that
$\lim_{k\rightarrow+\infty}\frac{P(\bar{x}+\tau_ku_k)-P(\bar x)}{\tau_k}=\xi$. Therefore we have
$$v=\lim_{k\rightarrow+\infty} v_k=-\lim_{k\rightarrow+\infty}\frac{P(\bar{x}+\tau_ku_k)-P(\bar{x})}{\tau_k} +\lim_{k\rightarrow+\infty}\frac{s_k-P(\bar{x})}{\tau_k}.$$
Since  $s_k \in \Lambda$, we have
$$\lim_{k\rightarrow+\infty}\frac{s_k-P(\bar{x})}{\tau_k}\in T_\Lambda(P(\bar{x})).$$
Hence $v\in -DP(\bar x)(u)+T_\Lambda(P(\bar{x}))$.

 $Case$ (ii) ($ \{ \frac{P(\bar{x}+\tau_ku_k)-P(\bar{x})}{\tau_k} \} $ is unbounded.) Without loss of generality, assume that  $$\lim_{k\rightarrow+\infty}\frac{\|P(\bar{x}+\tau_ku_k)-P(\bar x)\|}{\tau_k}=\infty.$$ Define
$t_k:=\|P(\bar{x}+\tau_ku_k)-P(\bar x)\|$.

Since $$\{ \frac{P(\bar{x}+\tau_ku_k)-P(\bar{x})}{t_k}\} =\{ \frac{P(\bar{x}+\tau_ku_k)-P(\bar{x})}{\|P(\bar{x}+\tau_ku_k)-P(\bar x)\|}\} $$ is bounded, we may  without loss of generality assume $\lim_{k\rightarrow+\infty}\{ \frac{P(\bar{x}+\tau_ku_k)-P(\bar{x})}{t_k}\}=\zeta.$
By definition of $DP(\bar{x})(0)$ and the fact that $\lim_{k\rightarrow \infty} \frac{\tau_k}{t_k}=0$, we have $$0\not = \zeta=\lim_{k\rightarrow+\infty}\frac{P(\bar{x}+\tau_ku_k)-P(\bar{x})}{t_k} =\lim_{k\rightarrow+\infty}\frac{P(\bar{x}+t_k (\frac{\tau_k}{t_k} u_k))-P(\bar{x})}{t_k} \in DP(\bar{x})(0).$$ Since $v_k\rightarrow v$ and $\lim_{k\rightarrow \infty} \frac{\tau_k}{t_k}=0$, we have
\begin{align*}
0=\lim_{k\rightarrow \infty}
\frac{\tau_k}{t_k} v_k&= -\lim_{k\rightarrow+\infty}\frac{P(\bar{x}+\tau_ku_k)-P(\bar{x})}{t_k} +\lim_{k\rightarrow+\infty}\frac{s_k-P(\bar{x})}{t_k}\\
&=-\zeta +\lim_{k\rightarrow+\infty}\frac{s_k-P(\bar{x})}{t_k} .
\end{align*}
Therefore $\zeta =\lim_{k\rightarrow+\infty}\frac{s_k-P(\bar{x})}{t_k}$ which implies that $\zeta \in T_\Lambda(P(\bar x))
$.

Conversely, let $v\in -DP(\bar{x})(u)+T_\Lambda(P(\bar{x}))$. Then  there exist $\xi\in DP(\bar{x})(u)$ and $\zeta\in T_\Lambda(P(\bar{x}))$ such that $v=-\xi+\zeta$.

If $P(x)$ is Hadamard directionally differentiable at $\bar x$,  then the limit
$$\xi=\lim_{t \downarrow 0, u'\rightarrow u}\frac{P(\bar x+t u')-P(\bar{x})}{t}$$
exists and  there exist sequences $\tau_k\downarrow 0,$ $s_k\xrightarrow{\Lambda}P(\bar{x})$ such that
$$\zeta=\lim_{k\rightarrow+\infty}\frac{s_k-P(\bar{x})}{\tau_k}.$$
Define  $$v_k=-\frac{P(\bar{x}+\tau_ku_k)-P(\bar{x})}{\tau_k}+\frac{s_k-P(\bar{x})}{\tau_k}=\frac{-P(\bar{x}+\tau_ku_k)+s_k}{\tau_k}.$$ {Then} $\lim_{k\rightarrow\infty}v_k=v$ and  $\tau_kv_k\in -P(\bar{x}+\tau_ku_k)+\Lambda$ for all $k$. Hence $v\in DG(\bar{x}|0)(u)$.

Now suppose that  $\Lambda$ is geometrically derivable. let $\tau_k\downarrow 0,\ u_k\rightarrow u$ be sequences such that
$$\xi=\lim_{k  \downarrow \infty }\frac{P(\bar x+\tau_k u_k)-P(\bar{x})}{\tau_k}.$$
Since  $\Lambda$ is geometrically derivable , there exists $s_k\in \Lambda$ such that
$$\zeta=\lim_{k\rightarrow+\infty}\frac{s_k-P(\bar{x})}{\tau_k}.$$
Define $$v_k=-\frac{P(\bar{x}+\tau_ku_k)-P(\bar{x})}{\tau_k}+\frac{s_k-P(\bar{x})}{\tau_k}=\frac{-P(\bar{x}+\tau_ku_k)+s_k}{\tau_k}.$$ {Then} $\lim_{k\rightarrow\infty}v_k=v$ and  $\tau_kv_k\in G(\bar{x}+\tau_ku_k)$ for all $k$. Hence  $v\in DG(\bar{x}|0)(u)$.
\endproof

\begin{defn}(coderivatives and directional coderivatives; see \cite[Definition 1.32]{Aub2} and \cite{BGO})
	For a set-valued map $G:\mathscr X\rightrightarrows\mathscr Y$ and a point $(\bar x, \bar y)\in gph G:=\{(x,y)\in\mathscr X\times\mathscr Y| y\in G(x)\}$,
  {the Fr\'{e}chet coderivative (Precoderivative) of $G$ at $(\bar x,\bar y)$ is a multifunction $\widehat D^*G(\bar x,\bar y):\mathscr
    Y \rightrightarrows\mathscr X$ defined as
    $$\widehat D^*G(\bar x,\bar y)(\zeta):=\left\{{\eta}\in\mathscr X\bigg|({\eta},-\zeta)\in \widehat N_{gph G}(\bar x,\bar y) \right\};$$}
	the limiting (Mordukhovich) coderivative of $G$ at $(\bar x,\bar y)$ is a multifunction $D^*G(\bar x,\bar y):\mathscr Y\rightrightarrows\mathscr X$ defined as
	$$D^*G(\bar x,\bar y)(\zeta):=\left\{{\eta}\in\mathscr X|({\eta},-\zeta)\in N_{gph G}(\bar x,\bar y) \right\}.$$
	The symbol $ D^*G(\bar x)$ is used when  $G$ is single valued.
	The limiting coderivative of $G$ at $(\bar x,\bar y)$ in direction $(u,\xi)\in\mathscr X\times\mathscr Y$ is defined as
	$$D^*G(\bar x,\bar y;(u,\xi))(\zeta):=\left\{{\eta}\in\mathscr X|({\eta},-\zeta)\in N_{gph G}(\bar x,\bar y;(u,\xi))\right\}.$$
	Similarly the symbol  $D^*G(\bar x;(u,\xi))$ is used when $G$ is single valued.
\end{defn}

\begin{remark} In the special case when $P:\mathscr{X}\rightarrow \mathscr{Y}$ is a single-valued map which is Lipschitz continuous at $\bar x$, by \cite[Theorem 3.28]{Aub2}, the coderivative is related to the limiting subdifferential in the following way:
$$ D^*P(\bar x)(\zeta)=\partial \langle P, \zeta\rangle (\bar x)\quad \mbox{for all}\ \zeta\in\mathscr{Y}.$$
By \cite[Proposition 5.1]{BGO}, if $P$ is Lipschitz near $\bar x$ in direction $u$, then $D^*P(\bar x;(u,\xi))(\zeta)\not=\emptyset$ if and only if $\xi\in DP(\bar x)(u)$, in which case
$$ D^*P(\bar x;(u,\xi))(\zeta)=\partial \langle P, \zeta\rangle (\bar x;(u,\langle\xi,\zeta\rangle)).$$

	Let $P:\mathscr X\rightarrow\mathscr Y$ be  $C^1$. By
	\cite[Remark 2.1]{BGO},  one has $DP(\bar x)(u)=\nabla P(\bar x)u$ and thus $D^*P(\bar x;(u,\xi))(\zeta)\not=\emptyset$ if and only if $\xi=\nabla P(\bar x)u$, in which case
	$$D^*P(\bar x;(u,\xi))(\zeta)=D^*P(\bar x)(\zeta)=\nabla P(\bar x)^*\zeta.$$
\end{remark}

To state our main results, given $P :\mathscr{X}\rightarrow \mathscr{Y}$ and $\Lambda \subseteq \mathscr{Y}$,  we define  the extended linearized cone as
\begin{equation} {\widetilde{\mathcal L}(x)}:=\left\{(u,\xi)\in\mathscr X\times\mathscr Y| \xi\in DP(x)(u)\cap T_\Lambda(P(x))\right\}.\label{Mset}
\end{equation} It is easy to see that the projection of $\widetilde{\mathcal L}(x)$ onto the space $\mathscr X$ is the linearized cone defined by
$\mathcal L(x):=\{u \in \mathscr X|\exists \xi \mbox{ such that } \xi\in DP(x)(u)\cap T_\Lambda(P(x))\}.$
When $P$ is differentiable at $x$, $DP(x)(u)=\nabla P(x)u$ and hence in this case
$$\widetilde{\mathcal L}(x)=\{(u, \nabla P(x)u):  0\in - \nabla P(x) u+T_\Lambda (P(x))\}.$$

\begin{prop}\label{Prop2.3}
	Let $P:\mathscr{X}\rightarrow \mathscr{Y}$ be continuous and $\Lambda \subseteq \mathscr{Y}$. Then
	\begin{equation}
	\widetilde{\mathcal L}(\bar x)=\{(0,0)\}\implies DG(\bar{x}|0)^{-1}(0)=\{0\}.
	\end{equation}
\end{prop}

\beginproof
By virtue of Proposition \ref{prop gd},  when $\widetilde{\mathcal L}(\bar x)=\{(0,0)\}$, one must have
$$ DG(\bar{x}|0)(u)\subseteq -DP(\bar x)(u)+T_\Lambda(P(\bar{x})).$$
Suppose that $u\in DG(\bar{x}|0)^{-1}(0)$.  Then equivalently, $0\in  DG(\bar{x}|0)(u)$.  Hence $0\in -DP(\bar x)(u)+T_\Lambda(P(\bar{x}))$ or equivalently $DP(x)(u)\cap T_\Lambda(P(x))\not =\emptyset$. Since $\widetilde{\mathcal L}(\bar x)=\{(0,0)\}$, it means that $\forall u\neq0,\ DP(x)(u)\cap T_\Lambda(P(x)) =\emptyset.$ Hence we must have $u=0$.
\endproof

\begin{prop}\label{strongMS} Let  $P:\mathscr X\rightarrow\mathscr Y$ be  continuous and  $\Lambda \subseteq  \mathscr Y$ be closed near $ \bar x\in \mathscr X$.  If $\widetilde{\mathcal L}(\bar x)=\{(0,0)\}$,  then $G(x)=P(x)-\Lambda$ is  strongly metrically subregular at $(\bar{x},0)$.
\end{prop}
\beginproof
 By \cite[Theorem 4C.1]{DonRock}, $G$ is strongly metrically subregular at $(\bar x, 0)$ if and only if $ DG(\bar{x}|0)^{-1}(0)=\{0\}$. The result then  follows from applying Proposition \ref{Prop2.3}.
\endproof

\section{{Weark} sufficient condition for metric subregularity}

In this section we will derive a  sufficient condition for metric subregularity of the system $P(x)\in \Lambda$ where $P(x):\mathscr X \rightarrow  \mathscr  Y$ is a continuous single-valued map and $\Lambda\subseteq \mathscr Y$ is locally closed.
Recall that no $\zeta$ satisfying condition $(\ref{con 1.1})$  alone is the so-called first-order sufficient condition for metric subregularity (FOSCMS) as established by Gfrerer and Klatte in  \cite[Corollary 1]{GKlatte16} for the case where $P$ is smooth and extended to the nonsmooth but calmness case in \cite[Proposition 2.2]{BGO}.   Our  sufficient condition in Theorem \ref{thm 2.1}  improves the FOSCMS in \cite[Proposition 2.2]{BGO} in two aspects.   First, we allow $P(x)$ to be only continuous instead of being calm. Secondly even in the case where $P(x)$ is calm, our condition is {weaker} in that the extra condition of the existence of sequences $ (u_k, v_k, \zeta_k)\rightarrow (u, 0, \zeta)\ \mbox{and}\ t_k\downarrow 0$  satisfying  (\ref{con 1.2}) and (\ref{con 1.3}) is required.

We will derive our result based on  the following  sufficient conditions for metric subregularity for general set-valued maps by Gfrerer  in \cite{Hcor}.

\begin{lemma}(see \cite[Corollary 1 and Remarks 1 and 2]{Hcor}\label{Hcor1})
	Let $G:\mathscr  X\rightrightarrows \mathscr Y$ be a closed set-valued map, and take a point $(\bar x,\bar y)\in gph G$.   Assume that  for any direction $u\in\mathscr X$, there do not exist sequences $t_k\downarrow 0, \|(u_k,v_k)\|=1, \|y_k^*\|=1$ with $\|u_k\|\rightarrow1, \|u\|u_k\rightarrow u, v_k\rightarrow0, x^*_k\rightarrow 0$ satisfying \[
	(x^*_k,-y^*_k)\in \widehat{N}_{gphG}(x_k',y_k'), \  x'_k\not \in G^{-1}(\bar y)
	\]
	and
	\[
	\lim_{k\rightarrow\infty}\frac{\langle y^*_k,y_k'-\bar y\rangle}{\|y'_k-\bar y\|}=1,
	\]
	where $x'_k:=\bar x+t_ku_k\neq\bar x,\ y'_k:=\bar y+t_kv_k\neq\bar y$. Then $G$ is metrically subregular  at $(\bar x,\bar y)$.
\end{lemma}
Note that as commented in \cite[Remark 2]{Hcor}, if the condition $x'_k\not \in G^{-1}(\bar y)$ is omitted then the resulting sufficient condition is stronger but may be easier to verify. However, in \cite[Example 1]{Ngai}, it was shown that sometimes these kinds of  conditions can not be omitted in order to show the metric subregularity.
\begin{lemma}\label{Hcod} Let  $P$ be a single-valued map from $\mathscr X$ to $\mathscr Y$ and $\Lambda$  be a subset of $\mathscr Y$. Define  $G(x):=P(x)-\Lambda,\  y= P( x)-s$ for some $ s\in \Lambda$. Then
$(x^*, -y^*)\in \widehat{N}_{gph G} ( x, y)$ implies that
$$ x^*\in \widehat{D}^* P( x)(y^*)  , \quad y^*\in \widehat{N}_\Lambda (P(x)-y).$$
\end{lemma}
\beginproof
Since $(x^*, -y^*)\in \widehat{N}_{gph G} ( x, y)$, by definition for any $\epsilon>0$,
\begin{equation}\label{ineq}
\langle x^*,x'-x \rangle+\langle -y^*,y'-y\rangle\leq\epsilon \|(x'-x,y'-y)\|
\end{equation}
 for any $(x',y')\in gph G$ which is sufficiently close to $(x,y)$. Let $y':= P(x)-s',\ s'\in \Lambda$. Then when $s'$ is close  to $s$, $y'=P(x)-s'$ is close to $y=P(x)-s$. Hence fixing  $x'=x$ in  $(\ref{ineq})$ we obtain that  for any $\epsilon>0$ and any $s'\in\Lambda$ sufficiently close to $s$,
\begin{equation*}
\langle -y^*,s-s'\rangle\leq\epsilon \|s-s'\|\Leftrightarrow\langle y^*,s'-s\rangle\leq\epsilon \|s-s'\|.
\end{equation*}
This means that $y^*\in \widehat N_\Lambda(s)=\widehat N_\Lambda(P(x)-y)$.

On the other hand, let $x'\in  \mathscr X$ and $y':=P(x')-s$. Then ${y'}\in G(x')$ and when $(x',P(x')) $ is close to $(x,P(x))$, $(x',y')$ is close to $(x,y)$.
Hence, by $(\ref{ineq})$ we have
\begin{equation*}
\langle x^*,x'-x \rangle+\langle -y^*,P(x')-P(x)\rangle\leq\epsilon \|(x'-x,P(x')-P(x))\|,
\end{equation*}
for any  $(x',P(x')) $  which is close to $(x,P(x))$.
This means that $$(x^*,-y^*)\in \widehat{N}_{gphP} ( x, P(x))$$ or equivalently $x^*\in \widehat{D}^* P( x)(y^*) $. The proof of the lemma is therefore complete.
\endproof

Applying Lemmas \ref{Hcor1} and \ref{Hcod}, we obtain the following sufficient condition for metric subregularity.
\begin{prop} \label{Prop3.1}Let  $P:\mathscr X\rightarrow\mathscr Y$ be a single-valued map and $\Lambda\subseteq \mathscr Y$ be closed. Let $G(x):=P(x)-\Lambda$ and   $ P(\bar x)\in \Lambda$.  Assume that $G(x)$ is a set-valued map which is closed around $\bar x$ and suppose that for any direction $u\in\mathscr X$, there do not exist sequences $t_k\downarrow 0,\ \|(u_k,v_k)\|=1,\ \|y_k^*\|=1$ with $\|u_k\|\rightarrow1,\ \|u\|u_k\rightarrow u,\ v_k\rightarrow0,\ x^*_k\rightarrow0$ satisfying \[
	x^*_k \in \widehat{D}^* P(\bar x+t_ku_k) (y^*_k), \quad y^*_k \in \widehat{N}_\Lambda (                                                                                                                                                                                                                                                                                                                   P(\bar x+t_ku_k)-t_kv_k), \quad { P(\bar x+t_ku_k)\not \in \Lambda}
	\]
	and
	\begin{equation}\label{seq1}
	 \lim_{k\rightarrow\infty}\frac{\langle y^*_k,v_k\rangle}{\|v_k\|}=1.
	\end{equation}
	Then $G$ is metrically subregular  at $(\bar x,0)$.
	
\end{prop}
Note that by \cite[Theorem 1.38]{Aub2}, when $P$ is Fr\'{e}chet differentiable  but not necessarily Lipschitz continuous, we have
$\widehat{D}^* P(x)(y^*)={\{\nabla P(x)^*y^*\}}.$

\begin{thm}\label{thm 2.1} Let $P:\mathscr X\rightarrow\mathscr Y$ be  continuous and $\Lambda \subseteq  \mathscr Y$ be closed at $ \bar x\in \mathscr X$. Suppose that  the weak sufficient condition for metric subregularity (WSCMS) holds at $\bar x$, i.e., for all $(0,0)\neq(u,\xi)\in \widetilde{\mathcal L}(\bar x)$, there exists no unit vector $\zeta$,  sequences $ (u_k, v_k, \zeta_k)\rightarrow (u, 0, \zeta)\ \mbox{and}\ t_k\downarrow 0$  satisfying
\begin{eqnarray}
&&0\in {D}^* P(\bar x; (u, \xi))(\zeta), \quad  \zeta\in N_\Lambda(P(\bar{x});\xi),\label{con 1.1}\\
&& 
\zeta_k\in \widehat{N}_\Lambda(s_k),\ s_k=P(\bar{x}+t_ku_k)-t_kv_k, \ { P(\bar x+t_ku_k)\not \in \Lambda},\label{con 1.2}\\
&&  \lim_{k\rightarrow \infty} \langle \zeta_k, \frac{v_k}{\|v_k\|}\rangle =1. \label{con 1.3}
\end{eqnarray}
Then  $G(x)=P(x)-\Lambda$ is metrically subregular at $(\bar{x},0)$.
\end{thm}

\beginproof
If $\widetilde{\mathcal L}(\bar x)=\{(0,0)\}$, then by Proposition \ref{strongMS}, $G$ is strongly metrically subregular  and hence metrically subregular at $(\bar x,0)$.
We now prove the result for the $\widetilde{\mathcal L}(\bar x)\not =\{(0,0)\}$ case by contradiction.
To the contrary, suppose that $P(x)-\Lambda$ is not metrically subregular at $(\bar x,0)$. By Proposition \ref{Prop3.1}, there exist $u\in\mathscr X$ and sequences $t_k\downarrow0, \|(u_k,v_k)\|=1,\ \|y_k^*\|=1$ with  $\|u_k\|\rightarrow1,\ \|u\|u_k\rightarrow u, v_k\rightarrow0,\ x^*_k\rightarrow0,$ such that
\begin{equation} \label{precode}
	(x^*_k, -y^*_k) \in \widehat{N}_{gph P}(\bar x+t_ku_k, P(\bar x+t_ku_k) ), y^*_k \in \widehat{N}_\Lambda (                                                                                                                                                                                                                                                                                                                   P(\bar x+t_ku_k)-t_kv_k)
	\end{equation}
	and  ${(\ref{seq1})}$ holds.
	
	Since we have $\|y_k^*\|=1,\ \|(u_k,v_k)\|=1$ and $v_k \rightarrow 0$, passing to a subsequence if necessary, we assume that $\lim_{k\rightarrow \infty} y_k^*=\zeta$, $\lim_{k\rightarrow \infty} u_k=u$ for certain $\|u\|=1$. It follows that $\|\zeta\|=1$.
	
	 $Case$ (1) ($\{ \frac{P(\bar{x}+t_ku_k)-P(\bar{x})}{t_k} \} $ is bounded.)
	Then without loss of generality we may assume that
$\lim_{k\rightarrow+\infty}\frac{P(\bar{x}+t_ku_k)-P(\bar x)}{t_k}=\xi$. Thus letting $\xi_k:=\frac{P(\bar{x}+t_ku_k)-P(\bar x)}{t_k}$, we have $P(\bar{x}+t_ku_k)=P(\bar x)+t_k\xi_k$. Combining with $(\ref{precode})$ we get
$$	(x^*_k, -y^*_k) \in \widehat{N}_{gph P}((\bar x,P(\bar x))+t_k(u_k,\xi_k)), y^*_k \in \widehat{N}_\Lambda (                                                                                                                                                                                                                                                                                                                   P(\bar x+t_ku_k)-t_kv_k).$$
Since $(u_k,\xi_k)\rightarrow(u,\xi)$ as $k\rightarrow\infty$, we have
$$(0,-\zeta) \in N_{gph P}((\bar x, P(\bar x)); (u,\xi)), \zeta \in N_\Lambda (P(\bar x); \xi).$$
  Also from the proof of Proposition \ref{prop gd},
  we see that $\xi\in DP(\bar x)(u)\cap T_\Lambda (P(\bar x))$ and hence $(u,\xi) \in \tilde {\mathcal L}(\bar x)$.

In summary for $Case$ (1), we have obtained a nonzero vector $\zeta$,  a nonzero vector $(u,\xi) \in \widetilde{L}(\bar x)$, and sequences  $(u_k,v_k,y^*_k)\rightarrow(u,0,\zeta)$ and $t_k\downarrow0$ such that
	\begin{align*}
	&0\in D^*P(\bar x;(u,\xi))(\zeta),\quad \zeta\in N_\Lambda(P(\bar x);\xi)\\
	&y^*_k\in\widehat N_{\Lambda}(s_k),\quad s_k=P(\bar x+t_ku_k)-t_kv_k,\\
	&\lim_{k\rightarrow\infty}\langle y^*_k,\frac{v_k}{\|v_k\|}\rangle=1,
	\end{align*}
	which contradicts the assumption in (WSCMS). Thus $P(x)-\Lambda$ is metrically subregular at $(\bar x,0)$.
	
 $Case$ (2) ($ \{ \frac{P(\bar{x}+t_ku_k)-P(\bar{x})}{t_k} \} $ is unbounded.)	 Without loss of generality, assume that  $\lim_{k\rightarrow+\infty}\frac{\|P(\bar{x}+t_ku_k)-P(\bar x)\|}{t_k}=\infty$.
	Define
	\begin{eqnarray*}
	&\tau_k:=\|t_ku_k\|+\|P(\bar x+t_ku_k)-P(\bar x)\|, &u'_k:=\frac{t_ku_k}{\tau_k},\\ &\xi_k:=\frac{P(\bar x+\tau_ku'_k)-P(\bar x)}{\tau_k}, &v'_k:=\frac{t_kv_k}{\tau_k}.
	\end{eqnarray*}

	Since $t_k/\tau_k \leq t_k/ \|P(\bar{x}+t_ku_k)-P(\bar{x})\|$,  we have $t_k/\tau_k\rightarrow 0$ and hence
	$v'_k\rightarrow0$ and 
	$u'_k\rightarrow0$. Since $\{\xi_k\}$ is bounded, taking a subsequence if necessary , we have
	\[
	 \xi:=\lim_{k\rightarrow\infty}\xi_k.
	\]	
	Then with $t_ku_k=\tau_ku'_k$ and $P(\bar x+t_ku_k)=P(\bar x)+\tau_k\xi_k$,
by $(\ref{precode})$  we get
	$$	(x^*_k, -y^*_k) \in \widehat{N}_{gph P}((\bar x,P(\bar x))+\tau_k(u'_k,\xi_k)), y^*_k \in \widehat{N}_\Lambda (                                                                                                                                                                                                                                                                                                                   P(\bar x+\tau_ku'_k)-\tau_kv'_k).$$
		
	Since $s_k=P(\bar x+\tau_ku'_k)-\tau_kv'_k$, we know that
	\begin{align*}
	\lim_{k\rightarrow\infty}\frac{s_k-P(\bar x)}{\tau_k}&=\lim_{k\rightarrow\infty}\frac{P(\bar x+\tau_ku'_k)-\tau_kv'_k-P(\bar x)}{\tau_k}\\
	&=\lim_{k\rightarrow\infty}\frac{P(\bar x+\tau_ku_k')-P(\bar x)}{\tau_k}=\xi.
	\end{align*}
	Thus, $\xi\in DP(\bar x)(0)\cap T_{\Lambda}(P(\bar x))$, which means $(0,\xi)\in \tilde {\mathcal L}(\bar x)$.
	With $x^*_k\rightarrow0$, we have
	\[
	0\in D^*P(\bar x;(0,\xi))(\zeta),\quad \zeta\in N_\Lambda(P(\bar x);\xi).
	\]
	
	By ${(\ref{seq1})}$, we can easily obtain that
	\[
	\lim_{k\rightarrow \infty} \langle y^*_k, \frac{v'_k}{\|v'_k\|}\rangle =\lim_{k\rightarrow \infty} \langle y^*_k, \frac{t_kv_k}{\|t_kv_k\|}\rangle=\lim_{k\rightarrow \infty} \langle y^*_k, \frac{v_k}{\|v_k\|}\rangle=1.
	\]
		
In summary for $Case$ (2), we obtain  a nonzero vector $\zeta$,  a nonzero vector $(0,\xi) \in \tilde {\mathcal L}(\bar x)$, and sequences  $(u'_k,v'_k,y^*_k)\rightarrow(0,0,\zeta)$ and $\tau_k\downarrow0$ such that
	\begin{align*}
	&0\in D^*P(\bar x;(0,\xi))(\zeta),\quad \zeta\in N_\Lambda(P(\bar x);\xi)\\
	&y^*_k\in\hat N_{\Lambda}(s_k),\quad s_k=P(\bar x+\tau_ku'_k)-\tau_kv'_k,\\
	&\lim_{k\rightarrow\infty}\langle y^*_k,\frac{v'_k}{\|v'_k\|}\rangle=1,
	\end{align*}
	which contradicts the assumption in (WSCMS). Thus $P(x)-\Lambda$ is metrically subregular at $(\bar x,0)$.
\endproof

{As an immediate consequence, if we discard the sequential conditions (\ref{con 1.2}) and (\ref{con 1.3}) in WSCMS, we derive from Theorem \ref{thm 2.1} the following sufficient condition for metric subregularity in the form of FOSCMS. The result improves \cite[Proposition 2.2]{BGO} in that $P$ is only assumed to be continuous instead of being  calm.}
\begin{cor} Let $P:\mathscr X\rightarrow\mathscr Y$ be  continuous and $\Lambda \subseteq  \mathscr Y$ be closed at $ \bar x\in \mathscr X$.
Suppose that  FOSCMS holds at $\bar x$, i.e.,
for all $ (u,\xi)$ such that $\xi\in DP(\bar x)(u)\cap T_\Lambda(P(\bar  x))$
$$0\in {D}^* P(\bar x; (u, \xi))(\zeta), \quad  \zeta\in N_\Lambda(P(\bar{x});\xi) \quad \Longrightarrow  \quad \zeta=0.$$
Then  $G(x)=P(x)-\Lambda$ is metrically subregular at $(\bar{x},0)$.
\end{cor}

\section{Directional quasi/pseudo normality}

As we mentioned in the introduction, quasi/pseudo-normality are also sufficient for metric subregularity. In this section we propose directional versions of the quasi-/pseudo-normality and show that they are slightly stronger than the WSCMS. Moreover we show that the SOSCMS implies pseudo-normality. Our results are based on the following observations.
\begin{prop}\label{cor 2.1} Let  $P:\mathscr X\rightarrow\mathscr Y$, $ (u^k, v^k, \zeta^k)\rightarrow (u, 0, \zeta), \ \ t_k\downarrow 0$ with $\|\zeta\|=1$, and
$s^k=P(\bar{x}+t_ku^k)-t_kv^k$.
	Then the condition
	\begin{equation}
	 \lim_{k\rightarrow \infty} \langle \zeta^k, \frac{v^k}{\|v^k\|}\rangle =1 \label{con 1.3new}
	 \end{equation}
	implies
	\begin{equation}
			 \zeta_i (P_i(\bar{x}+t_ku^k)-s_i^k)>0, \forall i\in I:=\{i: \zeta_i \not =0 \}
		 \label{con 1.5}
	\end{equation}
	which implies
	\begin{equation}
			 \langle \zeta,  P(\bar{x}+t_ku^k)-s^k \rangle>0.
			 \label{con 1.4}
	\end{equation}
	\end{prop}
	
\beginproof
 Suppose that $(\ref{con 1.3new})$ holds.
Since
\begin{eqnarray*}
\lefteqn{\left \|\frac{\zeta^k}{\|\zeta^k\|}-\frac{v^k}{\|v^k\|}\right \|^2}\\
&=& \langle \frac{\zeta^k}{\|\zeta^k\|}-\frac{v^k}{\|v^k\|}, \frac{\zeta^k}{\|\zeta^k\|}-\frac{v^k}{\|v^k\|} \rangle \\
&=&\frac{\|\zeta^k\|^2}{\|\zeta^k\|^2}-2\langle \frac{\zeta^k}{\|\zeta^k\|}, \frac{v^k}{\|v^k\|}\rangle+\frac{\|v^k\|^2}{\|v^k\|^2
}\\
&= & 2-\frac{2}{\|\zeta^k\|}\langle \zeta^k, \frac{v^k}{\|v^k\|}\rangle,
\end{eqnarray*}
$\lim_{k\rightarrow \infty} \langle \zeta^k, \frac{v^k}{\|v^k\|}\rangle =1$ and $\lim_{k\rightarrow\infty}\|\zeta^k\|=\|\zeta\|=1$, we have
 $$\lim_{k\rightarrow \infty}\|\frac{\zeta^k}{\|\zeta^k\|}-\frac{v^k}{\|v^k\|}\|=0.$$ Consequently, $\lim_{k\rightarrow \infty}\frac{v^k}{\|v^k\|}=\frac{\zeta}{\|\zeta\|}$. Thus when $k$ is large enough, for each $i=1,\dots,m$ with $\zeta_i\neq0$, $v^k_i$ has the same sign as $\zeta_i$. This means
$$\zeta_iv^k_i>0\quad  \forall i\in I:=\{i: \zeta_i\neq 0\},$$
which implies $(\ref{con 1.5})$.
Since $\zeta\neq0$,  $(\ref{con 1.5})$ obviously implies $(\ref{con 1.4})$.
\endproof

We are now in a position to define the concept of directional  quasi/pseudo-normality.
\begin{defn}[Directional  quasi/pseudo-normality] \label{quasinormal}Let   $P:\mathscr X\rightarrow\mathscr Y$ with $P(\bar x) \in \Lambda$.
\begin{itemize}
\item[(a)]We say that directional quasi-normality holds at $\bar{x}$ if  for all $$(0,0)\not = (u,\xi)\in \widetilde{\mathcal L}(\bar x):=\{(u,\xi)\in\mathscr X\times\mathscr Y| \xi\in DP(x)(u)\cap T_\Lambda(P(x))\},$$  there exists no $\zeta\neq 0$ such that
\begin{equation}\label{directionN}0\in D^*P(\bar x;(u,\xi))(\zeta),\quad   \zeta\in N_\Lambda(P(\bar{x});\xi)
	\end{equation}
	and
	\begin{eqnarray*}
	\left \{  \begin{array}{l}
	\exists (u^k,  s^k, \zeta^k)\rightarrow (u,  P(\bar{x}), \zeta)\mbox{and}\ t_k\downarrow0, \\
	\mbox{ s.t. } \zeta^k\in \widehat N_\Lambda(s^k)\mbox{ and }
	\zeta_i (P_i(\bar{x}+t_ku^k)-s^k_i)>0 \mbox{ if } \zeta_i \not =0.
	\end{array}\right.
	\end{eqnarray*}

\item[(b)] We say that directional pseudo-normality holds at $\bar{x}$ if   for all $(0,0)\not = (u,\xi)\in \widetilde{\mathcal L}(\bar x)$, there exists no $\zeta\neq0$ such that $(\ref{directionN})$ holds and
\begin{eqnarray*}
\left \{  \begin{array}{l}
\exists (u^k,  s^k, \zeta^k)\rightarrow (u,  P(\bar{x}), \zeta)\ \mbox{ and }\ t^k\downarrow0 ,\\
\mbox{ s.t. } \zeta^k\in \widehat N_\Lambda(s^k)\mbox{ and }
\langle\zeta, P(\bar{x}+t_ku^k)-s^k\rangle>0.
\end{array}\right.
\end{eqnarray*}
\end{itemize}
\end{defn}

{By virtue of Proposition \ref{cor 2.1}, directional pseudo-normality is stronger than directional quasi-normality. And consequently from Theorem \ref{thm 2.1}, they can provide sufficient conditions for metric subregularity.
}

\begin{cor}\label{thm 2.1new}
Let   $P:\mathscr X\rightarrow\mathscr Y$, $P(\bar x)\in \Lambda$, where  $P(x)$ is continuous at $\bar{x}$ and $\Lambda$ is closed near $\bar x$. If either directional pseudo-normality or directional quasi-normality holds at $\bar x$,
then  the set-valued map $G(x)=P(x)-\Lambda$ is metrically subregular at $(\bar{x},0)$.
\end{cor}
By definition, directional quasi-/pseudo-normality is weaker than quasi-/pseudo-normality, the following example shows that it is weaker than both quasi-normality and FOSCMS.

\begin{example}[FOSCMS fails but directional pseudo-normality holds]\label{eg4.1}
	Consider the constraint system defined by
$P(x)=(x,-x^2)\in \Lambda,$
	where $$\Lambda:=\{(x,y)|y\leq 0~or~y\leq x\}.$$
The point $\bar x=0$ is feasible since $(0,0)\in \Lambda$. We have
	\begin{align*}
	P(\bar{x})=
	(0,0),
	~~
	\nabla P(\bar{x})=
	\begin{pmatrix}
	1\\
-2\bar x\\
	\end{pmatrix}
	=
	\begin{pmatrix}
	1\\
	0
	\end{pmatrix}, T_\Lambda(P(\bar{x}))=\Lambda
	\end{align*}
	and the linearized cone $\mathcal L(\bar{x})=\{u\in\mathbb R| 0\in-\nabla P(\bar{x})u+T_\Lambda(P(\bar{x}))\}=\mathbb R$. Let $\bar u:=-1\in \mathcal L(\bar x)$, $\zeta:=(0,1)$ and $(x^k,y^k)=P(\bar{x})+\frac{1}{k}\nabla P(\bar{x})\bar u=(-\frac{1}{k},0)$. Then $\nabla P(\bar{x})^T\zeta=0$ and for each $k$, $\zeta\in N_\Lambda(x^k,y^k)$. Thus $\nabla P(\bar{x})^T\zeta=0$ and $\zeta\in N_\Lambda(P(\bar{x});\nabla P(\bar{x})\bar{u})$. Hence FOSCMS fails at $\bar{x}$.
	
	However, we can prove that directional pseudo-normality holds at $\bar x$. We prove it by contradiction. Assume that directional pseudo-normality fails at $\bar x$. Then there exist $0\not = u\in {\cal L}(\bar x)$,   $0\not= \zeta\in N_\Lambda(P(\bar x);\nabla P(\bar x)u)$ and a sequence $\{u^k,s^k,\zeta^k\}$ converging to $(u,P(\bar x),\zeta)$ and $t_k\downarrow 0$ such that
	\begin{equation}\label{eqn eg1}
	\nabla P(\bar x)^T\zeta=0,~\sum_{i=1}^2\zeta_i(P_i(\bar{x}+t_ku^k)-s^k_i)>0,~\zeta^k\in \widehat{N}_\Lambda(s^k).
	\end{equation}
	Solving $\nabla P(\bar{x})^T\zeta=0$, we obtain $\zeta_1=0$. Moreover since $N_\Lambda(P(\bar x))=\{0\}\times\mathbb R_+\cup\{(-r,r)|r\geq0\}$,  we have $\zeta\in \{0\}\times \mathbb R_{++}$. Since $\zeta^k\rightarrow \zeta$ and $\zeta^k \in \widehat{N}_\Lambda(s^k)$, we must have
	 $\zeta^k\in \{0\}\times \mathbb R_{++}$ and $s^k \in \{0\}\times \mathbb R_{+} $. Thus we obtain
	\begin{equation*}
	\sum_{i=1}^2\zeta_i(P_i(z^k)-s^k_i)
	=\zeta_2(P_2(z^k)-s^k_2)=\lambda(-(z^k)^2-s^k_2)\leq0,
	\end{equation*}
	where $z^k:=\bar{x}+t_ku^k$. But this
contradicts $(\ref{eqn eg1})$. Hence directional pseudo-normality holds at $\bar x$.
\end{example}

We now consider the case where  $\Lambda$ is the union of finitely many convex polyhedral sets in  $\mathscr{Y}$, i.e. $\Lambda:=\bigcup_{i=1}^p \Lambda_i$, where $$\Lambda_i:=\left \{y\in \mathscr{Y}|\langle \lambda_{ij},y \rangle\leq b_{ij},\quad j=1,...,m_i \right \}, \quad  i=1,...,p,$$ with  $\lambda_{ij}\in\mathscr{Y},\ b_{ij}\in\mathbb{R}$ for $j=1,...,m_i,$
are  convex polyhedral  sets. As noted in the introduction, by Robinson's multifunction theory \cite{Robinson81}, we know that when $P$ is affine and $\Lambda$ is the union of finitely many convex polyhedral sets, the set-valued map $G^{-1}$ is upper Lipschitz continuous and hence calm at each point of the graph. What is more, we now show that the { pseudo-normality always holds}. To our knowledge, this result has never been shown in the literature before.

The following results will be needed in the proof. For every $s\in \Lambda$, we denote by $p(s):=\{i=1,\ldots,p| s\in \Lambda_i\}$ the index set of the convex polyhedral sets containing $s$. Then we have from \cite{Gfr132} that
\begin{eqnarray}\label{tancon}
T_\Lambda(s)=\bigcup_{i\in p(s)} T_{\Lambda_i}(s),\ \quad \widehat{N}_\Lambda(s)=\bigcap_{i\in p(s)} \widehat{N}_{\Lambda_i}(s).
\end{eqnarray}

\begin{prop}\label{Prop linear}
	Let   $P:\mathscr X\rightarrow\mathscr Y$. Suppose that $P(x)$ is affine and $\Lambda$ is the union of finitely many convex polyhedral sets defined as above. Then for any feasible point $\bar{x}$ satisfying $P(\bar x)\in \Lambda$, pseudo-normality holds.
\end{prop}	
\beginproof
We prove the proposition by contradiction. 	{Assume that pseudo-normality does not hold at $\bar{x}$. Then there exists $ \zeta\neq0$ such that
\begin{eqnarray*}
	\left \{  \begin{array}{l}
		0=\nabla P(\bar{x})^*\zeta, \quad  \zeta\in N_\Lambda(P(\bar{x})),\\
		\exists (x^k, s^k, \zeta^k)\rightarrow (\bar x, P(\bar{x}), \zeta) \\
		\mbox{ s.t. } \zeta^k\in \widehat{N}_\Lambda(s^k),\
		\langle\zeta, P(x^k)-s^k\rangle>0.
	\end{array}\right.
\end{eqnarray*}}

As $s^k\rightarrow P(\bar{x})$ when $k\rightarrow\infty$ and $s^k\in \Lambda=\bigcup_{i=1}^p \Lambda_i$, by virtue of $(\ref{tancon})$, taking a subsequence if necessary, there exists $i\in \{1,\ldots, p\}$ such that for $k$ sufficiently large, $s^k\in \Lambda_i,\ P(\bar{x})\in \Lambda_i,\ \zeta^k\in N_{\Lambda_i}(s^k)$.   Define $J(s^k):=\{j=1,\ldots,m_i|\langle\lambda_{ij},s^k\rangle=b_{ij}\}$ and $J(P(\bar x)):=\{j=1,\ldots,m_i|\langle\lambda_{ij},P(\bar x)\rangle=b_{ij}\}$. Since $\zeta^k\neq0$, $s^k$ is not an interior point of $\Lambda_i$ and hence the index set $J(s^k)$ is not empty. Since $s^k\rightarrow P(\bar x)$, we have $J(s^k)\subseteq J(P(\bar x))$   when $k$ is sufficiently large. Hence without loss of generality, we can find a nonempty set $J\subseteq J(P(\bar x))$ such that $J(s^k)\equiv J$ for all $k$ large enough. Define $C:=\{\lambda_{ij}|j\in J\}$. Then we have $\zeta^k\in cone(C)$, where $$cone(C):=\{\Sigma_{j\in J}c_j\lambda_{ij}|c_j\geq0,\forall j\in J\}$$ denotes the conic hull of $C$. It follows that $\zeta\in cone(C)$. Since when $k$ large enough, for each $j\in J$, $\langle\lambda_{ij},P(\bar{x})-s^k\rangle=b_{ij}-b_{ij}=0$, we obtain
 $\langle\zeta, P(\bar{x})-s^k\rangle=0$.
Thus for sufficiently large $k$, we have
\begin{eqnarray*}
\lefteqn{0> \langle \zeta, s^k-P(x^k)\rangle+\langle \zeta, P(\bar{x})-s^k\rangle}\\
	&&=\langle \zeta, P(\bar{x})-P(x^k)\rangle\\
	&&=\langle \zeta, \nabla P(\bar{x})(\bar x-x^k)\rangle,
\end{eqnarray*}
which contradicts the condition that $ 0=\nabla P(\bar{x})^*\zeta$. Thus pseudo-normality holds at $\bar{x}$.
\endproof

For a single-valued mapping $P : \mathscr{X} \rightarrow \mathscr{Y}$ which is $C^1$ at $\bar x$ and $u\in \mathscr{Y}$, we define its second-order graphical derivative of $P(x)$ at $\bar x$  in direction $u$ as
\begin{eqnarray*}
&&D^2 P(\bar x)(u)\\
&&:=\left \{l\in\mathscr{Y}|\exists t_k\downarrow0,u^k\rightarrow u \ s.t.\ l=\lim_{k\rightarrow\infty}\frac{P(\bar x+t_ku^k)-P(\bar x)-t_k\nabla P({\bar x})u^k}{\frac{1}{2}t_k^2} \right \}.
\end{eqnarray*}

{In \cite[Theorem 4.3]{Gfr132},  a second-order sufficient condition for metric subregularity (SOSCMS) is presented for a split system in Banach spaces where one of the system is metrically subregular. Specializing  the result in \cite[Theorem 4.3]{Gfr132} to our system $(\ref{GS})$, we may conclude that if $P(x)$ is $C^1$ and  directionally second-order differentiable,  $\Lambda $ is   the union of finitely many convex polyhedral  sets and SOSCMS as stated in Theorem \ref{thm2.3} holds, then the system is directionally pseudo-normal.  In Theorem \ref{thm2.3},} {we extend this result to the case where $P(x)$ is $C^{1}$ and $\nabla P(x)$ is {directionally calm at $\bar x$ in each nonzero direction $u$ lying  in the linearization cone which means that  there exist positive numbers $\epsilon,\ \delta,\ L_u$ such that
	\begin{equation*}
	\|\nabla P(\bar x+tu')-\nabla P(\bar x)\|\leq L_u\|tu'\|\quad \mbox{for all}\ 0<t<\epsilon,\|u'-u\|<\delta.
	\end{equation*}}
Moreover we show that SOSCMS implies directional pseudo-normality.

\begin{thm}\label{thm2.3}
	{Let $P(\bar x )\in \Lambda$ where $P(x)$ is $C^{1}$, $\Lambda $ is the union of finitely many convex polyhedral  sets in $\mathscr{Y}$ and $\nabla P(x)$ is directionally calm at $\bar x$ in each direction $0\not =u$
	such that $\nabla P(\bar x)u \in T_\Lambda(P(\bar x))$.}  Suppose  SOSCMS holds at $\bar x$, i.e.,   for all $0\neq u$ such that $\nabla P(\bar x)u \in T_\Lambda(P(\bar x))$,    there exists no $\zeta\neq0$ such that 
	\begin{eqnarray*}
		\left \{  \begin{array}{l}
		 \nabla P(\bar{x})^* \zeta =0, \ \zeta\in N_\Lambda(P(\bar{x});\nabla P(\bar{x})u),\\
			\exists l\in D^2P(\bar{x})(u)\ s.t.\ \langle \zeta,l\rangle>0.\\
		\end{array}\right.
	\end{eqnarray*}
	Then $\bar x$ is directionally pseudo-normal  at $\bar x$.
\end{thm}
\beginproof
We prove that SOSCMS is stronger than directional pseudo-normality by contradiction. Assume there exist  $0\neq u$ such that $\nabla P(\bar x)u \in T_\Lambda(P(\bar x))$ and  $\zeta\neq0$ such that
\begin{eqnarray*}
\left \{  \begin{array}{l}
 \nabla P(\bar{x})^* \zeta=0, \quad  \zeta\in N_\Lambda(P(\bar{x});\nabla P(\bar{x})u)\\
\exists (u^k,  s^k, \zeta^k)\rightarrow (u,  P(\bar{x}), \zeta)~\mbox{and}~ t_k\downarrow0 \\
\mbox{ s.t. } \zeta^k\in \widehat{N}_\Lambda(s^k),\
\sum_{i=1}^m\zeta_i (P_i(\bar{x}+t_ku^k)-s^k_i)>0.
\end{array}\right.
\end{eqnarray*}

 Notice that $\langle P(x),e_j\rangle$, where $e_j$ is in  the orthogonal basis $\mathscr{E}$,  is a function on $\mathscr X$. By the mean value theorem, for each $j$ and $k$, there exist $0<c^k_j<t_k$ such that
$$\langle P(\bar{x}+t_ku^k)-P(\bar{x}), e_j\rangle=\langle\nabla P(\bar{x}+c_j^ku^k)(\bar{x}+t_ku^k-\bar{x}),e_j\rangle=\langle\nabla P(\bar{x}+c_j^ku^k)t_ku^k,e_j\rangle.$$
Thus
\begin{eqnarray*}
	\lefteqn{\left\langle\frac{P(\bar{x}+t_ku^k)-P(\bar{x})-t_k\nabla P(\bar{x})u^k}{\frac{1}{2}t_k^2},e_j\right\rangle}\\
	&&=\frac{1}{2t_k}\left(\frac{\langle P(\bar{x}+t_ku^k)-P(\bar{x}),e_j\rangle}{t_k}-\langle\nabla P(\bar{x})u^k,e_j\rangle\right)\\
	&&=\frac{2}{t_k}\left(\frac{\langle \nabla P(\bar{x}+c^k_ju^k)u^k,e_j\rangle t_k}{t_k}-\langle\nabla P(\bar{x})u^k,e_j\rangle\right)\\
	&&=\frac{2}{t_k}(\langle \nabla P(\bar{x}+c^k_ju^k)u^k,e_j\rangle-\langle\nabla P(\bar{x})u^k,e_j\rangle).
\end{eqnarray*}
{Since $\nabla P(x)$ is directionally calm at $\bar x$ in direction $u$, there exists  $L_u> 0 $ such that  for each $j$ and sufficiently large $k$,}
\begin{eqnarray*}
\lefteqn             {\left \|\frac{2}{t_k}(\langle \nabla P(\bar{x}+c^k_ju^k)u^k,e_j\rangle-\langle\nabla P(\bar{x})u^k,e_j\rangle)\right \|}\\
&&\leq \frac{2L_u\|\bar{x}+c^k_ju^k-\bar{x}\|\|u^k\|}{t_k}\\
&&\leq\frac{2L_ut_k\|u^k\|^2}{t_k}=2L_u\|u^k\|^2.
\end{eqnarray*}
This implies that the sequence $\left \{\langle\frac{P(\bar{x}+t_ku^k)-P(\bar{x})-t_k\nabla P(\bar{x})u^k}{\frac{1}{2}t_k^2},e_j\rangle \right \}$ is bounded. Consequently, the sequence $\left \{\frac{P(\bar{x}+t_ku^k)-P(\bar{x})-t_k\nabla P(\bar{x})u^k}{\frac{1}{2}t_k^2} \right \}$ is bounded. Taking a subsequence if necessary, there exists $l$ such that
$$l:=\lim_{k\rightarrow\infty}\frac{P(\bar{x}+t_ku^k)-P(\bar{x})-t_k\nabla P(\bar{x})u^k}{\frac{1}{2}t_k^2}\in  D^2P(\bar{x})(u).$$
It follows that
\begin{eqnarray}
0&&<\langle \zeta, P(\bar{x}+t_ku^k)-s^k\rangle
\nonumber\\
&&=\langle \zeta, P(\bar{x}+t_ku^k)-P(\bar{x})+P(\bar{x})-s^k\rangle \nonumber \\
&&=\langle \zeta,t_k\nabla P(\bar{x})u^k+\frac{t_k^2}{2}l+o(t_k^2)\rangle+\langle \zeta, P(\bar{x})-s^k\rangle.\label{eqn comp}
\end{eqnarray}

By assumption, $\nabla P(\bar{x})^*\zeta=0$, which means $\langle \zeta,t_k\nabla P(\bar{x})u^k\rangle=0$. And since $s^k\rightarrow P(\bar{x})$ as $k\rightarrow\infty$, taking a subsequence if necessary, there exists $j\in \{1,\dots, p\}$ such that for $k$ sufficiently large, $s^k\in \Lambda_j,\ P(\bar{x})\in \Lambda_j,\ \zeta^k\in N_{\Lambda_j}(s^k)$. Since $\Lambda_j$ is convex polyhedral, similar to the discussion in the proof of Proposition $\ref{Prop linear}$, we have $\langle \zeta, P(\bar{x})-s^k\rangle=0$. Thus for $k$ large enough, by $(\ref{eqn comp})$ we have
\begin{eqnarray*}
	0&&<\langle \zeta,t_k\nabla P(\bar{x})u^k+\frac{t_k^2}{2}l+o(t_k^2)\rangle+\langle \zeta, P(\bar{x})-s^k\rangle\\
	&&\leq \frac{t_k^2}{2}\langle \zeta,l+o(1)\rangle.
\end{eqnarray*}
Then we obtain that $\exists l\in D^2P(\bar{x})(u)$ such that $\langle \zeta,l \rangle\geq0$. But this contradicts the assumption of the SOSCMS. The contradiction proves that the  SOSCMS implies directional pseudo-normality.
\endproof

{Since directional calmness is obviously weaker than calmness, the following corollary follows immediately from Theorem \ref{thm2.3}. We  say that $P(x)$ is $C^{1,c}$ at $\bar x$ if $P(x)$ is $C^1$ at $\bar x$ and $\nabla P(x)$ is calm at $\bar x$, i.e., there exist $\kappa >0$ and a neighborhood $U$ of $\bar x$ such that $\|\nabla P(x)-\nabla P(\bar x)\|\leq \kappa \|x-\bar x\|$ for all $x\in U$.}
{\begin{cor}Let $P(\bar x )\in \Lambda$ where $P$ is $C^{1,c}$ and $\Lambda $ is the union of finitely many convex polyhedral sets in $\mathscr{Y}$. Suppose  SOSCMS holds at $\bar x$. Then $\bar x$ is directionally pseudo-normal.
\end{cor}}

In summary,  we have shown the following implications:
\begin{eqnarray*}
&&  SOSCMS\implies \mbox{directional\ pseudo-normality} \implies \mbox{directional quasi-normality} \\
&& \implies \mbox{WSCMS} \implies \mbox{metric subregularity}.
\end{eqnarray*}

\section{Applications to  complementarity and KKT systems}
In this section we apply our results to complementarity and KKT systems. When directional quasi-/pseudo-normality are applied to these systems we derive expressions that are much simpler and moreover can be directly compared with classical  quasi/pseudo-normality.

First we consider the complementarity system formulated as follows:
\begin{eqnarray*}
\mbox{(CS)}\ && H(x)=0,  \quad 0 \leq \Phi(x)\perp \Psi(x)\geq 0,
\end{eqnarray*}
where $ H(x): \mathbb{R}^n \rightarrow \mathbb{R}^d,\ \Phi,\ \Psi: \mathbb{R}^n \rightarrow \mathbb{R}^m$.
For simplicity of explanation, we omit possible  inequality and abstract constraints and moreover we assume that
all functions are continuously differentiable. The results can be extended to the general case in a straightforward manner.

Define $\Omega_{EC}:=\{(a,b)\in\mathbb R_+\times\mathbb R_+|  ab=0\}$. For any set $C$ and any positive integer $m$ we denote by $C^m$ the $m$-Cartesian product of $C$.  (CS) can be rewritten as
\begin{eqnarray*}
	&&\left ( H(x), (\Phi_1(x), \Psi_1(x)),\ldots, (\Phi_m(x), \Psi_m(x) )\right )
	\in  \{0\}^d \times \Omega_{EC}^m.
\end{eqnarray*}
To derive the precise form of the directional quasi-/pseudo-normality, we review the  formulas for the regular normal cone, the limiting normal cone, the tangent cone and the directional limiting normal cone of the set $\Omega_{EC}$.
\begin{lemma}\label{VG}\cite[Lemma 4.1]{Gfr14}
	The Fr\'{e}chet normal cone  to $\Omega_{EC}$ is
	\begin{eqnarray*}
		\widehat N_{\Omega_{EC}}(a,b)=
		\left\{
		-(\gamma, \nu)\bigg|\begin{array}{ll}
			\nu=0,\quad &\mbox{if}\ 0=a<b\\
			\gamma \geq0,\nu \geq0,\quad &\mbox{if}\ a=b=0\\
			\gamma=0,\quad &\mbox{if}\ a>b=0
		\end{array}
		\right\},
	\end{eqnarray*}
	the limiting normal cone is
	\begin{eqnarray*}
		N_{\Omega_{EC}}(a,b)=
		\left\{
		\begin{array}{ll}
			\widehat N_{\Omega_{EC}}(a,b)\quad &\mbox{if}\ (a,b)\neq(0,0)\\
			\{-(\gamma,\nu)|\mbox{either }\gamma>0,\nu>0\ \mbox{or }\ \gamma\nu=0\},\quad &\mbox{if}\ (a,b)=(0,0)
		\end{array}
		\right\},
	\end{eqnarray*}
	and the tangent cone is
	\begin{eqnarray*}
		T_{\Omega_{EC}}(a,b)=
		\left\{
		(d_1,d_2)\bigg|\begin{array}{ll}
		d_1=0,\quad &\mbox{if}\ 0=a<b\\
			(d_1,d_2)\in\Omega_{EC},\quad &\mbox{if}\ a=b=0\\
			d_2=0,\quad &\mbox{if}\ a>b=0
		\end{array}
		\right\}.
	\end{eqnarray*}
	
	For all $d=(d_1,d_2)\in T_{\Omega_{EC}}(a,b)$, the directional limiting normal cone to $\Omega_{EC}$ in direction $d$ is
	\begin{eqnarray*}
		N_{\Omega_{EC}}((a,b);d)=
		\left\{
		\begin{array}{ll}
			N_{\Omega_{EC}}(a,b)\quad &\mbox{if}\ (a,b)\neq(0,0)\\
			N_{\Omega_{EC}}(d_1,d_2),\quad &\mbox{if}\ (a,b)=(0,0)
		\end{array}
		\right\}.
	\end{eqnarray*}
\end{lemma}

Let $\bar x$ be a feasible point of the system (CS). We define index sets
\begin{align*}
I_{00}:=I_{00}(\bar x)&:=\{i|\Phi_i(\bar x)=0,\Psi_i(\bar x)=0\},\\
I_{0+}:=I_{0+}(\bar x)&:=\{i|\Phi_i(\bar x)=0,\Psi_i(\bar x)>0\},\\
I_{+0}:=I_{+0}(\bar x)&:=\{i|\Phi_i(\bar x)>0,\Psi_i(\bar x)=0\},
\end{align*}
and define the linearized cone as
\begin{eqnarray*}
	{\mathcal L}(\bar x):=
	\left\{
	u \in\mathbb R^n|\begin{array}{ll}
		0=\nabla H_i(\bar x)u  & i=1,\dots, d, \\
		0=\nabla \Phi_i(\bar x)u &  i\in I_{0+}, \\
		 0=\nabla \Psi_i(\bar x)u\quad & i\in I_{+0},\\
		(\nabla \Phi_i(\bar x)u,\nabla \Psi_i(\bar x)u)\in\Omega_{EC},\quad &i\in I_{00}
	\end{array}
	\right\}.
\end{eqnarray*}

Given $u\in \mathcal L(\bar x)$ we define
\begin{align*}
& I_{+0}(u):=\{i\in I_{00}| \nabla {\Phi_i}(\bar x)u>0=\nabla  {\Psi}_i(\bar x)u\},\\
& I_{0+}(u):=\{i\in I_{00}| \nabla \Phi_i(\bar x)u=0<\nabla \Psi_i(\bar x)u\},\\
& I_{00}(u):=\{i\in I_{00}|\nabla \Phi_i(\bar x)u=0=\nabla \Psi_i(\bar x)u \}.
\end{align*}
Let $\bar x$ be a feasible point of (CS).  By Definition \ref{quasinormal} and Proposition \ref{prop YZ},  since the complementarity set $\Omega_{EC}$ is directionally regular,  (CS) is  directionally quasi- or pseudo-normal if and only if for all directions $0\neq u\in \mathcal L(\bar x)$ there exists no $(\eta, \gamma, \nu )\neq 0$ such that
	\begin{eqnarray}
	 &&  0=\nabla H(\bar x)^T \eta-\nabla \Phi(\bar x)^T \gamma -\nabla \Psi(\bar x)^T \nu, \label{CSquasinormal1}\\
		&&  -(\gamma_i, \nu_i) \in N_{\Omega_{EC}}(\Phi_i(\bar x), \Psi_i(\bar x);\nabla \Phi_i(\bar x) u, \nabla \Psi_i(\bar x) u) \quad i=1,\dots, m, \qquad  \label{CSquasinormal2}\\
	&&   \exists (u^k, h^k,\phi^k, \psi^k, \eta^k,\gamma^k,\nu^k)\rightarrow (u, H(\bar x),  \Phi(\bar{x}),\Psi(\bar x), \eta, \gamma,\nu), t_k\downarrow0 \nonumber \\
	&& \mbox{  such that }
		\left\{\begin{array}{ll}
		\eta^k \in N_{\{0\}^d}(h^k),{-}(\gamma_i^k, \nu_i^k) \in \widehat{N}_{\Omega_{EC}}(\phi_i^k,\psi_i^k)\quad i=1,\dots, m,\\
		\mbox{if}\ \eta_i\neq0,\ \eta_i(H_i(\bar x+t_ku^k)-h_i^k)>0,\\
		\mbox{if}\ \gamma_i\neq0,\ \gamma_i(\Phi_i(\bar x+t_ku^k)-\phi_i^k)<0,\\
		\mbox{if}\ \nu_i\neq0,\ \nu_i(\Psi_i(\bar x+t_ku^k)-\psi_i^k)<0.
		\end{array}\right. \label{CSquasinormal3}
		\end{eqnarray}
		or $$ \eta^T(H(\bar x+t_ku^k) -h^k)-\gamma^T(\Phi (\bar x+t_ku^k)-\phi^k)-\nu^T (\Psi(\bar x+t_ku^k)-\psi^k)>0  $$ respectively.
	
By the formula for the directional limiting normal cone in Lemma \ref{VG}, $(\ref{CSquasinormal2})$ is equivalent to (ii) in the following definition. Since $\eta^k \in N_{\{ 0\}^d}({h^k})$, we have ${h^k=0}$.  Suppose $\gamma_i\not =0$. Then for sufficiently large $k$, $\gamma_i^k\not =0$. Since  ${-}(\gamma_i^k,\ \nu_i^k) \in \widehat{N}_{\Omega_{EC}}(\phi_i^k,\psi_i^k)$ we must have  $\phi_i^k=0$. Similarly if $\nu_i\not =0$, we must have $\psi_i^k=0$.
Based on these discussions, the directional quasi-normality for (CS) can be written in the following form which is much more concise.
\begin{defn}\label{QuasiCS}
	Let $\bar x$ be a feasible solution of (CS). $\bar x$ is said to be directionally quasi- or pseudo-normal if for all directions $0\neq u\in \mathcal L(\bar x)$ there exists no $(\eta,\gamma, \nu )\neq 0$ such that
	\begin{itemize}
		\item[(i)] $ 0=\nabla {H}(\bar x)^T \eta-\nabla {\Phi}(\bar x)^T \gamma -\nabla {\Psi}(\bar x)^T \nu;$\\
		\item[(ii)]
	$\gamma_{i}=0,\ i\in I_{+0}\cup I_{+0}(u)$;\ 	$\nu_{i}=0,\ i\in I_{0+}\cup I_{0+}(u)$;\ eihter $\gamma_i>,\ \nu_i> 0$ or $\gamma_i \nu_i=0$, $ i\in I_{00}(u)$;\\
		\item[(iii)]   $\exists u^k\rightarrow u$ and $t_k\downarrow0$ such that
		\begin{eqnarray*}
		\left\{\begin{array}{ll}
		\mbox{if}\ \eta_i\neq0,\ \eta_iH_i(\bar x+t_ku^k)>0\\
		\mbox{if}\ \gamma_i\neq0,\ \gamma_i\Phi_i(\bar x+t_ku^k)<0,\\
		\mbox{if}\ \nu_i\neq0,\ \nu_i\Psi_i(\bar x+t_ku^k)<0.
		\end{array}\right.
		\end{eqnarray*}
	or $$ \eta^TH(\bar x+t_ku^k) -\gamma^T\Phi (\bar x+t_ku^k)-\nu^T \Psi(\bar x+t_ku^k)>0  $$ respectively.
	\end{itemize}
\end{defn}
\begin{remark} In Definition \ref{QuasiCS}, if we only require that there exists no $(\eta,\gamma,\nu)\neq0$ satisfying condition $(i)$,  then it reduces to the linearly independent constraint qualification (MPEC-LICQ) {(see \cite{SS})}. If we only require that there exists no $({\eta},\gamma,\nu)\neq0$ satisfying condition $(i)$ and  change (ii) to $$\gamma_{i}=0,\ i\in I_{+0};\ 	
\nu_{i}=0,\ i\in I_{0+},\
 \mbox{ eihter } \gamma_i\geq0,\ \nu_i\geq0\ \mbox{  or }\ \gamma_i \nu_i=0,\ i\in I_{00}$$ then it reduces to MPEC-NNAMCQ {\cite[Definition 2.10]{Ye2005}}. If we omit (iii), then it reduces to FOSCMS. If we take $u$ to be any direction,  then it reduces to the MPEC quasi-/pseudo-normality first  given in   \cite[Definition 3.2]{Kanzow} and extended to the Lipschitz continuous case in \cite[Definition 5]{YZ2014}.
  Since for the set $\Omega_{EC}$ and any $0\neq d\in T_{\Omega_{EC}}(0,0)$, the directional normal cone $N_{\Omega_{EC}}((0,0);d)$ is strictly smaller than the limiting normal cone $N_{\Omega_{EC}}(0,0)$, if there exists some $u\in \mathcal L(\bar x)$ such that $(\nabla G(\bar x)u,\nabla H(\bar x)u)\neq(0,0)$, then directional quasi-/pseudo-normality will be strictly weaker than standard quasi-/pseudo-normality.
\end{remark}

We now consider the following KKT system of an optimization problem with equality and inequality constraints:
\begin{eqnarray*}
&& \nabla_x L(x,\mu, \lambda)=0,\\
&& \mu\geq 0,\quad g(x)\leq 0,\quad \langle g(x),\mu\rangle=0,\\
&& h(x)=0,
\end{eqnarray*}
where $f:\mathbb{R}^p\rightarrow \mathbb{R},\ g:\mathbb{R}^p\rightarrow \mathbb{R}^m,\ h:\mathbb{R}^p\rightarrow \mathbb{R}^n$ are twice continuously differentiable,  $\mu\in \mathbb{R}^m,\ \lambda\in \mathbb{R}^n$, and $L(x,\mu,\lambda):=f(x)+\mu^Tg(x)+\lambda^Th(x)$ is the Lagrange function. Denote  the feasible set of the KKT system by $\mathcal{F}_{KKT}$.
We say that the error bound property holds at $(x^*,\mu^*,\lambda^*)\in \mathcal{F}_{KKT}$ if
 there exist $\alpha>0$ and $U$, a neighborhood {of  $(x^*,\mu^*,\lambda^*)$}, such that
\begin{eqnarray}\label{serrorb}
d_{{\cal F}_{KKT}}(x,\mu,\lambda) \leq \alpha \max\{ \|\nabla_x L(x,\mu,\lambda)\|,\|h(x)\|, \|\min \{\mu, -g(x)\}\| \}, \quad \forall (x,\mu,\lambda)\in U.
\end{eqnarray}
It is easy to see that this error bound property can be derived from the metric subregularity/calmness of KKT system and hence directional quasi-/pseudo-normality is a sufficient condition.
Such an error bound property
is crucial to the quadratic convergence of the Newton-type method (see \cite{FFH}).  The classical sufficient conditions for the error bound property are either  MFCQ combined with  the second-order sufficient condition (SOSC) or requiring $g,\ h $ to be affine and $f$ to be quadratic  (see e.g., \cite{r82}). These sufficient conditions were weakened in \cite{f02,hg} but still require SOSC.  Recently weaker sufficient conditions have been proposed  including the existence of noncritical multipliers, a concept introduced by Izmailov for pure equality constraint in \cite{i05}, extended by Izmailov and Solodov \cite[Defnition 2]{is} to problems with inequalities and further extended to a general variational system by Mordukhovich and Sarabi \cite[Definition 3.1]{ms}.
Note that  as shown in \cite[Proposition 3]{is}, the existence of noncritical multipliers is equivalent to a stronger type of error bound property: existence of $\alpha>0$ and $U$, a neighborhood of  $(x^*,\mu^*,\lambda^*)$, such that
\begin{eqnarray*}
\|x-\bar x\|+d_{{\cal M}(\bar{x})}(\mu,\lambda) \leq \alpha \max\{ \|\nabla_x L(x,\mu,\lambda)\|,\|h(x)\|, \|\min \{\mu, -g(x)\}\| \} \quad \forall (x,\mu,\lambda)\in U,
\end{eqnarray*}
where ${\cal M}(\bar x):=\{(\mu,\lambda): 0=\nabla_x L(\bar x, \mu, \lambda),\  \mu\geq 0,\ \langle g(\bar x),\
\mu\rangle=0\}$ denotes the set of multipliers. Obviously this is a stronger error bound property
  than the error bound property $(\ref{serrorb})$.

The KKT system is a special case of (CS) with
\[
H(x,\mu,\lambda):=(\nabla_x L(x,\mu,\lambda)
, h(x)), \quad \Phi(x,\mu,\lambda):=-g(x),\quad  \Psi(x,\mu,\lambda):=\mu.
\]
Let $(\bar x,\bar \mu,\bar \lambda)$ be a feasible point of the KKT system. We define the following index sets:
\begin{align*}
I_{00}:=I_{00}(\bar x,\bar \mu,\bar \lambda)&:=\{i|g_i(\bar x)=0,\ \bar \mu_i=0\},\\
I_{0+}:=I_{0+}(\bar x,\bar \mu,\bar \lambda)&:=\{i|g_i(\bar x)=0,\ \bar \mu_i>0\},\\
I_{+0}:=I_{+0}(\bar x,\bar \mu,\bar \lambda)&:=\{i|-g_i(\bar x)>0,\ \bar \mu_i=0\}.
\end{align*}
The linearization cone for the KKT system is
\begin{eqnarray*}
	\mathcal L(\bar x,\bar \mu,\bar \lambda):=
	\left\{
	u=(u^x,u^\mu,u^\lambda) |\begin{array}{ll}
		0=\nabla_{xx}^2 L(\bar x,\bar u,\bar v)u^x+\nabla g(\bar x)^Tu^\mu+\nabla h(\bar x)^T u^\lambda &\\
		0=\nabla h(\bar x)u^x &\\
		0=\nabla g_i(\bar x)u^x,\quad &i\in I_{0+}\\
		0= u^{\mu}_i,\quad &i\in I_{+0}\\
		u^\mu_i\geq0, \nabla g_i(\bar x)u^x\leq0\ \mbox{and}\ u^\mu_i\nabla g_i(\bar x)u^x=0,\quad &i\in I_{00}
	\end{array}
	\right\}.
\end{eqnarray*}
Given $u\in {\mathcal L}(\bar x,\bar \mu,\bar \lambda)$ we define the index sets
\begin{align*}
& I_{+0}(u):=\{i\in I_{00}| -\nabla g_i(\bar x)u^x>0=u^\mu_i\},\\
& I_{0+}(u):=\{i\in I_{00}| \nabla g_i(\bar x)u^x=0<u^\mu_i\},\\
& I_{00}(u):=\{i\in I_{00}| \nabla g_i(\bar x)u^x=0=u^\mu_i \}.
\end{align*}
Then by Definition \ref{QuasiCS}, we propose the following definition of directional quasi-normality for the KKT system.
\begin{defn}\label{QuasiKKT}
Let $(\bar x,\bar \mu,\bar \lambda)$ be a feasible point of the KKT system. $(\bar x,\bar \mu,\bar \lambda)$ is said to be directionally quasi-/pseudo-normal if for all directions $$0\neq \bar{u}:=(\bar u^x,\bar u^{\mu},\bar u^\lambda)\in {\mathcal L}(\bar x,\bar \mu,\bar \lambda)$$ there exists no $(\xi, \zeta, \eta)\neq 0$ such that
\begin{itemize}
	\item[(i)] $ 0=\nabla^2_{xx}L(\bar x,\bar \mu,\bar \lambda)\xi+\nabla h(\bar x)^T\eta+\nabla g(\bar x)^T\zeta$;\\
	\item[(ii)] $\nabla h(\bar x)\xi=0$;
	\item[(iii)]
		$\zeta_{i}=0,\ i\in I_{+0}\cup I_{+0}(\bar u)$;\ 	$\nabla g_{i}(\bar x)\xi=0,\ i\in I_{0+}\cup I_{0+}(\bar u)$;\ either $\zeta_i>0,\ \nabla g_i(\bar x)\xi>0$ or $\zeta_i \nabla g_i(\bar x)\xi=0$, $ i\in I_{00}(\bar u)$;\\
		\item[(iv)]   $\exists u_k:=(u_k^x,u_k^{\mu},u_k^\lambda)\rightarrow \bar u$ and $t_k\downarrow0$ such that
		\begin{eqnarray*}
			\left\{\begin{array}{ll}
				\mbox{if}\ \xi_i\neq0,\ \xi_i\nabla_x L_i((\bar x,\bar \mu,\bar \lambda)+t_ku_k)>0,\\
				\mbox{if}\ \eta_i\neq0,\ \eta_ih_i(\bar x+t_ku^x_k)>0,\\
				\mbox{if}\ \zeta_i\neq0,\ \zeta_ig_i(\bar x+t_ku^x_k)>0,\\
				\mbox{if}\ (\nabla g(\bar x)\xi)_i\neq0,\ (\nabla g(\bar x)\xi)_i(\bar u^\mu_i+t_k(u^\mu_k)_i)<0,
			\end{array}\right.
		\end{eqnarray*}	
		or $$\xi^T \nabla_x L((\bar x,\bar \mu, \bar \lambda)+t_ku_k) +\eta^T h (\bar x+t_k  u_k^x) -( \nabla g(\bar x)\xi)^T (\bar u+ t_k u_k^\mu) >0 $$ respectively.
	\end{itemize}
\end{defn}

\begin{remark} Let $(\bar x,\bar \mu,\bar \lambda)$ be a feasible point to the KKT system.
By  \cite[Definition 2]{is}, $(\bar \mu,\bar \lambda)\in {\cal M}(\bar x)$   is a critical multiplier associated with $\bar x$ if there exists $(\xi, \zeta, \eta)$ with $\xi\neq0$ satisfying that
\begin{eqnarray*}
	\left\{\begin{array}{ll}
		 0=\nabla^2_{xx}L(\bar x,\bar \mu,\bar \lambda)\xi+\nabla h(\bar x)^T\eta
		+\nabla g(\bar x)^T\zeta,&\\
	 0=\nabla h(\bar x)\xi,\\
		0=\nabla g_i(\bar x)\xi,\quad &i\in I_{0+}\\
		0= \zeta_i,\quad &i\in I_{+0}\\
		\zeta_i\geq0, \nabla g_i(\bar x)\xi\leq0\ \mbox{and}\ \zeta_i\nabla g_i(\bar x)\xi=0,\quad &i\in I_{00}.
		\end{array}
	\right.
\end{eqnarray*}
Note that from Definition \ref{QuasiKKT}, we can see that  even if $(\bar \mu,\bar \lambda)$ is a  critical multiplier with $\bar x$, it is still possible for directional quasi-normality to hold. In particular let  $(\xi, \zeta, \eta)$ satisfy Definition \ref{QuasiKKT} with $\xi\not =0$. Suppose that
for $ i\in I_{00}(\bar u)$, it is not possible to have $\zeta_i>0,\ \nabla g_i(\bar x)\xi>0$. Then $(\bar \mu,\bar \lambda)\in {\cal M}(\bar x)$   is a critical multiplier associated with $\bar x$.
\end{remark}

\section*{Acknowlegement} The authors would like to thank Helmut Gfrerer for an inspiring discussion on an earlier version of this paper, {and the anonymous referees for their helpful suggestions and comments.}


\end{document}